\documentclass[11pt]{article}

\usepackage{amsmath,amsthm}

\usepackage{amssymb,latexsym}

\usepackage{enumerate}

\newtheorem{tw}{Theorem}[section]
\newtheorem{lm}[tw]{Lemma}
\newtheorem{wn}[tw]{Corollary}
\newtheorem{stw}[tw]{Proposition}
\newenvironment{dow}{\it Proof.\rm}{\hfill $\Box$}

\theoremstyle{definition}
\newtheorem{df}[tw]{Definition}
\newtheorem{uw}[tw]{Remark}
\newtheorem{prz}[tw]{Example}

\newcommand{\BN}{{\mathbb N}}
\newcommand{\BR}{{\mathbb R}}
\newcommand{\BX}{{\mathbb X}}

\newcommand{\FF}{{\mathcal{F}}}

\newcommand{\MM}{{\mathcal{M}}}

\newcommand{\EE}{{\mathcal{E}}}

\newcommand{\esssup}{\mathop{\mathrm{ess\,sup}}}

\numberwithin{equation}{section}

\begin{document}
\title{Systems of quasi-variational inequalities related
to the switching problem}
\author{Tomasz Klimsiak}
\date{}
\maketitle

\begin{abstract}
We prove  the existence  of weak solution for a system of
quasi-variational inequalities related to a switching problem with
dynamic driven by operator associated with a semi-Dirichlet form
and with measure data. We give a stochastic representation of
solutions in terms of solutions of a system of reflected BSDEs
with oblique reflection. As a by-product, we prove the existence
of an optimal strategy in the switching problem and show
regularity of the payoff function.
\end{abstract}
\noindent {\small\bf Mathematics Subject Classification (2010):}
35J75, 60J45.

\footnotetext{T. Klimsiak: Institute of Mathematics, Polish
Academy of Sciences, \'Sniadeckich 8, 00-956 Warszawa, Poland, and
Faculty of Mathematics and Computer Science, Nicolaus Copernicus
University, Chopina 12/18, 87-100 Toru\'n, Poland. E-mail:
tomas@mat.umk.pl}

\section{Introduction}
\label{sec1}

Let $E$ be a locally compact separable metric space, $m$ be a
Radon measure on  $E$ with full support, and let $(L, D(L))$ be
the generator of  a regular semi-Dirichlet form $(\EE,D[\EE])$ on
$L^2(E;m)$. The class of such operators is quite wide. The model
example of local operator associated with semi-Dirichlet form is
the second order uniformly elliptic divergence form operator with
bounded drift, i.e. operator of the form
\begin{equation}
\label{eq0.01} L=\sum_{i,j=1}^d \frac{\partial}{\partial x_j}
\Big(a_{ij}(\cdot)\frac{\partial}{\partial x_i}\Big)
+\sum_{i=1}^db^i(\cdot)\frac{\partial}{\partial x_i}\,.
\end{equation}
As the example of nonlocal operator of this class may serve
\begin{equation}
\label{eq1.6} L=\Delta^{\alpha(\cdot)},
\end{equation}
i.e.  fractional Laplacian with possibly varying exponent
$\alpha:E\rightarrow(0,2)$ satisfying some regularity assumptions.

In the paper we consider the following problem: for given
functions $f^j:E\times \BR^N\rightarrow\BR$,
$h_{j,i}:E\times\BR\rightarrow\BR$, $i,j=1,\dots,N$,  smooth (with
respect to the capacity associated with $(\EE,D[\EE])$) measures
$\mu^1,\dots,\mu^N$ on $E$ and sets $A_1,\dots,A_N$ such that
$A_j\subset \{1,\dots, j-1,j+1,\dots,N\}$ find a pair $(u,\nu)$
consisting of a function $u=(u^1,\dots,u^N):E\rightarrow\BR^N$ and
a vector $\nu=(\nu^1,\dots,\nu^N)$ of smooth measures on $E$ such
that
\begin{equation}
\label{eq1.1} \left\{
\begin{array}{l}-Lu^{j}
=f^{j}(\cdot,u)+\mu^j+\nu^{j},\quad \medskip\\
\int_E(u^{j}-\max_{i\in A_j}h_{j,i}(\cdot,u^{i}))\,d\nu^{j}=0,\medskip\\
\, u^{j}\ge\max_{i\in A_j}h_{j,i}(\cdot,u^{i}),\quad j=1,\dots,N.
\end{array}
\right.
\end{equation}
Intuitively, we are looking for $u$ satisfying the equations
$-Lu^{j} =f^{j}(\cdot,u)+\mu^j$ on the sets $\{u^{j}>\max_{i\in
A_j}h_{j,i}(\cdot,u^{i})\}$. The measure $\nu^j$ represents the amount
of energy we have to add to the system to keep $u^j$ above the
obstacle $H^j(\cdot,u):=\max_{i\in A_j}h_{j,i}(\cdot,u^{i})$. The second
equation in (\ref{eq1.1}) says that $\nu^j$ is minimal in the
sense that it acts only when $u^j=H^j(\cdot,u)$.

Systems of the form (\ref{eq1.1}) arise  when considering the
so-called switching problem. They were studied by many authors
(see, e.g., \cite{DHP,EH,HJ,HM,HT,HZ,HZhao,LNO,LNO1}) in  case $L$ is a
diffusion operator or  diffusion operator perturbed by nonlocal
operator associated with a Poisson measure, and the data are
$L^2$-integrable  (hence, in particular, $\mu^i=0$,
$i=1,\dots,N$). Also note that in all the papers cited above $f$
is Lipschitz continuous with respect to $u$ and viscosity
solutions are considered.

In the present paper we generalize the existing results on
(\ref{eq1.1}) in the sense that we consider quite general class of
operators and  measure data. We also considerably weaken
the assumptions on $f$, because we only assume that it is
quasi-monotone with respect to $u$.

When $h_{j,i}$ do not depend on $u$, system (\ref{eq1.1})
resembles the usual system of variational inequalities written in
complementary form (see \cite{KS} and also
\cite{K:PA,K:arx,KR:JFA} for the case of one equation). Such a
form has proved useful in the study of variational inequalities
with measure data (see \cite{K:arx,KR:JEE,RS}). One of the main
reason is that it allows one to use known  results on semilinear
elliptic PDEs with measure data. On the other hand, the usual
variational approach is applicable only to systems with
$L^2$-data.

Our general approach to (\ref{eq1.1}) (system of quasi-variational
inequalities in complementary form) is similar to that in
\cite{K:arx,KR:JEE}. It can be briefly described as follows. Let
$\BX=(\{X_t,\,t\ge0\},\{P_x,x\in E\})$ be a Hunt process with life
time $\zeta$ associated with $(\EE,D[\EE])$, and for smooth
measure $\gamma$ let $A^\gamma$ denote the continuous additive
functional of $\BX$ in the Revuz correspondence with $\gamma$. By
a solution of (\ref{eq1.1}) we mean a pair $(u,\nu)$ satisfying
the second and third condition in (\ref{eq1.1}), and such that for
quasi-every $x\in E$ (with respect to the capacity associated with
$\EE$) the following generalized nonlinear Feynman-Kac formula is
satisfied
\begin{equation}
\label{eq1.5} u(x)=E_x\Big(\int_0^\zeta
f(X_r,u(X_r))\,dr+\int_0^\zeta \,dA^\nu_r
+\int_0^\zeta\,dA^\mu_r\Big).
\end{equation}
Note that from (\ref{eq1.1}) one can often deduce some regularity
properties of $u$. For instance, if $\mu$ is a measure of bounded
variation, $(u,\nu)$ satisfies (\ref{eq1.5}) and we know that
$f(\cdot,u)\in L^1(E;m)$ and $\nu$ has also bounded variation,
then $T_ku\in D_e[\EE]$ for every $k>0$, where $D_e[\EE]$ is the
extended Dirichlet space for $\EE$  and $T_ku(x)=((-k)\vee
u(x))\wedge k$. In fact, $(u,\nu)$ is then a renormalized solution
of the first equation in (\ref{eq1.1}) in the sense introduced in
\cite{KR:NDEA} (for the case where $L$ is of the form
(\ref{eq0.01}) see also \cite{DMOP} and \cite{RS}).

Roughly speaking, to find a solution $(u,\nu)$ of (\ref{eq1.1}) in
the sense described above  we find a solution of some system of
Markov-type BSDEs with oblique reflection associated with
(\ref{eq1.1}), and we study various properties of these solutions.
Then, using some ideas from the papers \cite{KR:JFA,KR:CM} devoted
to PDEs with measure data, we translate the results on this
systems of reflected BSDEs into results on (\ref{eq1.1}).

As a matter of fact, in the first part of the paper we study
general, nonMarkov-type BSDEs. First, in Section \ref{sec2}, we
give an existence result for solutions of system of BSDEs of the
type
\[
Y^{j}_{t}=\xi^{j}+\int_{t}^{T}f^{j}(r,Y_r)\,dr +\int_t^T
dV^j_r-\int_{t}^{T}\,dM_r^j,\quad t\in [0,T],
\]
where $V$ is a finite variation c\`adl\`ag process, with
quasi-monotone  right-hand side $f$, i.e. off-diagonal increasing
and on-diagonal decreasing. This type of equation was not
considered in the literature in such generality. Then, in Section
\ref{sec3}, we prove the existence of a solution of the system of
RBSDEs with oblique reflection of the form
\begin{equation}
\label{eq1.3} \left\{
\begin{array}{l}
Y^{j}_{t}=\xi^{j}+\int_{t}^{T}f^{j}(r,Y_r)\,dr +\int_t^T
dV^j_r+\int_t^T dK^{j}_r-\int_{t}^{T}\,dM_r^j,\quad t\in [0,T] ,
\medskip\\
Y^{j}_{t}\ge \max_{i\in A_j}h_{ji}(t,Y^{i}_{t}),\quad t\in [0,T],
\medskip\\
\int_0^T (Y^{j}_{t}-\max_{i\in A_j}
h_{j,i}(t,Y^{i}_{t}))\,dK^{j}_{t}=0, \,j=1,\dots,N.
\end{array}
\right.
\end{equation}
This result generalizes the  existence results proved for
$L^2$-data and Brownian filtration (see, e.g., \cite{HZ})  or
filtration generated by a Brownian motion and an independent
Poisson measure (see \cite{HZhao,LNO}) to the case of general
filtration and $L^1$-data. Moreover, as compared with
\cite{HZ,HZhao,LNO}, we impose less restrictive assumptions on the
off-diagonal growth of the right-hand side. We also allow the
terminal time $T$ to be unbounded stopping time. In Section
\ref{sec3} we also show that solution of (\ref{eq1.3}) may be
approximated by solutions of the  system of penalized BSDEs
\begin{align*}
Y^{n,j}_t=\xi^j+\int_t^T f^j(r,Y^n_r)\,dr+\int_t^T\,dV^j_r
+\int_t^Tn(Y^{n,j}_r-H^j(r,Y^n_r))^-\,dr-\int_t^T\,dM^{n,j}_r
\end{align*}
with $H^j$ of the form
\[
H^j(t,y)= \max_{i\in A_j}h_{j,i}(t,y^i).
\]

In Section \ref{sec4} we study  the switching problem (we describe
it briefly below) and its connection with reflected BSDEs.
Therefore we restrict our attention to $h_{j,i}$ of the form
\begin{equation}
\label{eq1.8} h_{j,i}(t,y)=c_{j,i}(t)-y^i
\end{equation}
for some adapted continuous processes $c_{j,i}$ (in applications
$c_{j,i}(t)$ is  the cost of switching the process of, say
production,  from mode $j$ to mode $i$ in time $t$). Our main
result says that  if $f$ in (\ref{eq1.3}) does not depend on $y$
then the first component $Y$ of the solution of (\ref{eq1.3})  is
the value function of the switching problem.

In Section \ref{sec5}, using the results of the probabilistic part
of the paper, we first give an existence result for (\ref{eq1.1}),
and we show that $u$ may be approximated by solutions of the
following system of penalized PDEs
\[
-Lu_n^j=f^j(\cdot,u_n)+n(u_n^j-H^j(\cdot,u_n))^{-}+\mu
\]
with
\[
H^j(x,y)= \max_{i\in A_j}h_{j,i}(x,y^i).
\]
We also give conditions ensuring that $f(\cdot,u)\in L^1(E;m)$ and
the measures $\nu^j$ have bounded variation. In particular, under
these conditions,  $T_k(u^j)\in D_e[\EE]$ and  $u^j$ is a
renormalized solution of the first equation in (\ref{eq1.1}) (see
comment following (\ref{eq1.5})).  We next turn to the switching
problem of Section \ref{sec4}, but in the Markovian setting, i.e.
in case $f^j(t,y)=f^j(X_t)$, $c_{j,i}(t)=c_{j,i}(X_t)$ for some
$f^j,c_{j,i}:E\rightarrow\BR$. The problem can be described as
follows. Consider a factory  in which we can change a mode of
production. Let $c_{j,i}(X)$ be the cost of the change from  mode
$j$ to mode $i$, and let $\psi_i(X)+dA^{\mu^i}$ be the payoff rate
in mode $i$. Then a  management strategy
$\mathcal{S}=(\{\tau_n\},\, \{\xi_n\})$ consist of a pair of two
sequences of random variables. The variable $\tau_n$ is the moment
when we decide to switch the mode of production, and $\xi_n$ is
the mode to which we switch at time $\tau_n$. If $\xi_0=j$ then we
start the production at  mode $j$. Under strategy $\mathcal{S}$
the expected profit on the interval $[0,T]$ is given by the
formula
\[
J(x,\mathcal{S},j)=E_x\Big(\int_0^T \psi_{w_r}(X_r)\,dr
+\int_0^T\,dA^{\mu^{w_r}}_r-\sum_{n\ge 1} c_{w_{\tau_{n-1}},
w_{\tau_{n}}}(X_{\tau_n})\mathbf{1}_{\{\tau_n<\zeta\}}+\xi^{w_T}\Big),
\]
where
\[
w_t=\xi_0\mathbf{1}_{[0,\tau_1)}(t)+\sum_{n\ge 1} \xi_n
\mathbf{1}_{[\tau_{n},\tau_{n+1})}(t).
\]
A strategy $\mathcal{S}^*$ is called optimal if
$J(x,\mathcal{S}^*,j)=\sup_{\mathcal{S}}J(x,\mathcal{S},j)$. In
Section \ref{sec5} we show that under some assumption on the data
there exists an optimal strategy, and moreover, if $T=\zeta$, then
$u$ defined by the formula
\[
u^j(x)=J(x,\mathcal{S}^*,j)
\]
is a unique solution of (\ref{eq1.1}) with $h_{j,i}$ defined by
(\ref{eq1.8}).

\section{Systems of BSDEs with quasi-monotone generator}
\label{sec2}

Let $(\Omega,\FF,P)$ be a probability space, $\mathbb{F}=
\{\FF_t,\, t\ge 0\}$ be a filtration satisfying the usual
conditions,  and let $T$ be a stopping time.

 In what follows $N\in\BN$,
$\xi=(\xi^1,\dots,\xi^N)$ is  an $\FF_T$-measurable random vector,
$V=(V^1,\dots,V^N)$ is an $\mathbb{F}$-adapted process such that
$V_0=0$  and each component $V^j$ is a process of finite
variation, $f:\Omega\times[0,T]\times\BR^N\rightarrow\BR^N$ is a
measurable function such that for every $y\in\BR^N$ the process
$f(\cdot,y)$ is $\mathbb{F}$-progressively measurable. As usual,
in the sequel, in our notation we will omit the dependence of $f$
on $\omega\in\Omega$.

We set $|V|_t=\sum^N_{j=1}|V^j|_t$, where $|V^j|_t$ stands for
that variation of $V^j$ on $[0,t]$, and we adopt the following
notation:
\[
f^j(t,y;a)=f^j(t,y_1,\dots,y_{j-1},a,y_{j+1},\dots,y_N),
\quad y\in\BR^N,\, a\in\BR
\]
and
\[
\underline{f}^j(t,a)=\inf_{y\in\BR^N}f^j(t,y;a),\quad
\overline{f}^j(t,a)=\sup_{y\in\BR^N}f^j(t,y;a),\quad a\in\BR.
\]

For $x=(x^1,\dots,x^N)$ we set $|x|=\sum_{j=1}^N |x^j|$, and for
$x,y\in\BR^N$ we write $x\le y$ if $x^j\le y^j$, $j=1,\dots N$.
For  processes $X,Y$ we write $X\le Y$ if $X_t\le Y_t$, $t\in
[0,T\wedge a]$ for all $a\ge 0$, and $X=Y$ if $X\le Y$ and $X\ge
Y$. The abbreviation ucp means``uniformly on compacts in
probability".

The following assumptions will be needed throughout the paper.
\begin{enumerate}
\item[(A1)] $E(|\xi|+\int_0^T\,d|V|_r)<\infty$,
\item[(A2)]  for every $t\in [0,T]$ the function $f(t,\cdot)$ is on-diagonal
decreasing, i.e. for $j=1,\dots,N$ we have $f^j(t,y;a)\le
f^j(t,y;a')$ for all $ a\ge a',\, a,a'\in\BR$, $y\in\BR^N$,
\item[(A3)] for every $t\in [0,T]$ the function $f(t,\cdot)$ is
off-diagonal increasing,
i.e. for $j=1,\dots,N$ we have $f^j(t,y)\le f^j(t,y')$ for all
$y,y'\in\BR^N$ such that $y\le y'$ and $y^j=y'^{j}$,
\item[(A4)] $y\mapsto f(t,y)$ is continuous for every $t\in [0,T]$,
\item[(A5)] $\int_0^{T\wedge a}|f(r,y)|\,dr<\infty$ for all $y\in\BR^N,\, a\in\BR$.
\end{enumerate}

Note that functions satisfying (A2) and (A3) are called
quasi-monotone.

Recall that an adapted  c\`adl\`ag process $\eta$ is of class (D)
if the collection $\{\eta_{\tau}: \tau$ is a finite valued stoping
time\} is uniformly integrable.

\begin{df}
We say that a pair $(Y,M)$ of $N$-dimensional $\mathbb{F}$-adapted
processes is a solution of the system of backward stochastic
differential equations on the interval $[0,T]$ with terminal
condition $\xi$ and right-hand side $f+dV$ (BSDE$^T(\xi,f+dV)$ for
short) if
\begin{enumerate}
\item[(i)] $Y^j$ is of class (D), $M^j$ is a local martingale such that
$M^j_0=0$, $j=1,\dots,N$,
\item[(ii)] $\int_0^{T\wedge a}|f(r,Y_r)|\,dr<\infty$ for every $a\ge 0$,
\item[(iii)] for $j=1,\dots, N$ and all $a\ge 0$,
\[
Y^j_t=Y^j_{T\wedge a}+\int_t^{T\wedge a}
f^j(r,Y_r)\,dr+\int_t^{T\wedge a} dV^j_r-\int_t^{T\wedge a}\,
dM^j_r,\quad t\in [0,T\wedge a].
\]
\item[(iv)] $Y_{T\wedge a}\rightarrow \xi$ $P$-a.s.
as $a\rightarrow \infty$.
\end{enumerate}
\end{df}

\begin{uw}
\label{uw3.0} Let $(Y,M)$ be a solution of BSDE$^T(\xi,f+dV)$. If
\begin{equation}
\label{eq2.01}
E\Big(|\xi|+\int_0^T|f(r,Y_r)|\,dr+\int_0^T\,d|V|_r\Big)<\infty
\end{equation}
then $M$ is a uniformly integrable martingale and
\begin{equation}
\label{eq2.1}
Y^j_t=E(\xi^j+\int_t^Tf^j(r,Y_r)\,dr+\int_t^T\,dV^j_r|\FF_t),\quad
t\le T,\quad j=1,\dots,N.
\end{equation}
To see this, we set $\tilde M=(\tilde M^1,\dots,\tilde M^N)$,
where
\[
\tilde{M}^j_t=E\Big(\int_0^T f^j(r,Y_r)\,dr
+\int_0^T\,dV^j_r|\FF_t\Big)-Y^j_0.
\]
An elementary computation shows that $(Y,\tilde{M})$ is a solution
of BSDE$^T(\xi,f+dV)$. Hence $M=\tilde{M}$. Therefore we may pass
to the limit as $a\rightarrow\infty$ in condition (iii) of the
above definition. We then get
\[
Y^j_t=\xi^j+\int_t^T f^j(r,Y_r)\,dr+\int_t^T
dV^j_r-\int_t^TdM^j_r, \quad t\le T,\quad j=1,\dots,N.
\]
Since $M$ is a uniformly integrable martingale, this yields
(\ref{eq2.1}).
\end{uw}

\subsection{One-dimensional equations}
\label{sec2.1}

In this subsection we assume that $N=1$.
\begin{uw}
\label{uw3.123} Let $\eta_t=E(\xi|\FF_t)$,
$f_\eta(t,y)=f(t,y+\eta_t)$. If a pair $(\bar{Y},\bar{M})$ is a
solution of BSDE$^T(0,f_\eta+dV)$, then the pair $(Y,M)$ defined
by
\[
Y_t=\bar{Y}_t+\eta_t,\quad M_t=\bar{M}_t+\eta_t-\eta_0
\]
is a solution of BSDE$^T(\xi,f+dV)$.
\end{uw}

\begin{stw}
\label{tw3.123} Let $\eta_t=E(\xi|\FF_t)$,
$t\ge0$. If {{\rm(A1), (A2), (A4), (A5)}} are satisfied, and
moreover,
\begin{equation}
\label{eq3.123} E\int_0^T|f(r,\eta_r)|\,dr<\infty,
\end{equation}
then there exists a solution of {\rm BSDE}$^T(\xi,f+dV)$.
\end{stw}
\begin{dow}
Let $f_\eta(t,y)=f(t,y+\eta_t)$. Then by \cite[Theorem
3.4]{KR:JFA} there exists a solution  $(\bar{Y},\bar{M})$ of
BSDE$^T(0,f_\eta+dV)$, and hence,  by Remark \ref{uw3.123}, there
exists a solution of BSDE$^T(\xi,f+dV)$.
\end{dow}

Assumption (\ref{eq3.123}) is quite natural  in the theory of
BSDEs with random terminal time (see, e.g., \cite{BDHPS}). We
would like, however, to weaken it and show that in fact
assumptions {{\rm(A1), (A2), (A4), (A5)}} together with
(\ref{eq3.123}) holding true with  some semimartingale $\eta$ of
class (D) and integrable finite variation part are sufficient for
the existence of a solution of BSDE$^T(\xi,f+dV)$. That
(\ref{eq3.123}) can be weaken is quite easy to see in case $T$ is
finite. In the general case more work have to be done.

\begin{uw}
Condition (\ref{eq3.123}) is too strong in many important
application. To illustrate, let us consider the well known
penalization scheme for reflected BSDE with terminal condition
$\xi=0$, coefficient equal to zero and lower barrier $L$, that is
equation of the form
\begin{equation}
\label{eq2.02} Y^n_t=\int_t^Tn(Y^n_r-L_r)^-\,dr-\int_t^T\,dM^n_r.
\end{equation}
Of course, this is  BSDE$^T(0,f_n)$ with $f_n(t,y)=n(y-L_t)^-$.
Suppose that $L_t=t^{-1}\mathbf{1}_{[1,\infty)}(t)$ and
$T=\infty$. We then expect that there exists a solution
$(Y^n,M^n)$ of (\ref{eq2.02}) and $\{Y^n\}$ converges to the Snell
envelope of $L$, which exists since $L$ is of class (D). Observe
that in this example  (\ref{eq3.123}) does not hold with
$\eta_t=E(\xi|\FF_t)=0$. However, (\ref{eq3.123}) is satisfied
with $\eta$ replaced by the semimartingale $L$. The same
phenomenon can happen for finite $T$. To see this, let us consider
a finite stoping time  $\tau$ such that $\tau\ge 1$ and $E\ln
\tau=\infty$, and set $T=\tau+1$. Let
$L_t=t^{-1}\mathbf{1}_{[1,\tau)}(t)$. Then (\ref{eq3.123}) is not
satisfied with $\eta=0$, but is satisfied with $\eta$ replaced by
the semimartingale $L$.
\end{uw}

\begin{lm}
\label{lm3.123} If {{\rm (A2), (A4), (A5)}} are
satisfied and
\begin{equation}
\label{eq2.3} E\big(|\xi|+\int_0^T
d|V|_r+\int_0^T|f(r,0)|\,dr\big)^2<\infty,
\end{equation}
then there exists a solution of {\rm BSDE}$^T(\xi,f+dV)$.
\end{lm}
\begin{dow}
Let $g$ be a strictly positive function on $\BR^+$ such that
$\int_0^\infty g(r)\,dr<\infty$. Write
\[
f_{n,m}=(f\wedge (n \cdot g))\vee(-m\cdot g).
\]
By Proposition \ref{tw3.123}, for all $n,m\in\mathbb{N}$ there is
a solution $(Y^{n,m},M^{n,m})$ of BSDE$^T(\xi,f_{n,m}+dV)$. By
\cite[Proposition 3.1]{KR:JFA}, $Y^{n,m}\le Y^{n+1,m}$. Set
$Y^m_t=\sup_{n\ge 1} Y^{n,m}_t$. Applying the Tanaka-Meyer formula
and (A2) we get
\[
|Y^{n,m}_t|\le E\big(|\xi|+\int_0^T|f(r,0)|\,dr +\int_0^T\,
d|V|_r|\FF_t\big)=:X_t,\quad t\le T.
\]
By (\ref{eq2.3}),  $E\sup_{t\ge 0}|X_t|^2<\infty$, whereas by
\cite[Lemma 2.3, Lemma 2.5]{KR:JFA},
\begin{equation}
\label{eq3.1234}
\sup_{n,m\ge 1} E\big(\int_0^T|f_{n,m}(r,Y^{n,m}_r)|\,dr\big)^2<\infty.
\end{equation}
By Remark \ref{uw3.0} we have
\[
Y^{n,m}_t=E\big(\xi+\int_t^Tf(r,Y^{n,m}_r)\,dr+\int_t^T\,dV_r|\FF_t\big),\quad
t\le T.
\]
Letting $n\rightarrow\infty$ in the above inequality and using
(\ref{eq3.1234}) we obtain
\begin{equation}
\label{eq2.4}
Y^{m}_t=E\big(\xi+\int_t^Tf(r,Y^{m}_r)\,dr+\int_t^T\,dV_r|\FF_t\big).
\end{equation}
Set
\[
M^m_t=E\big(\xi+\int_0^Tf(r,Y^{m}_r)\,dr+\int_0^T\,dV_r|\FF_t\big)-Y_0^m,
\quad t\le T.
\]
Then the pair $(Y^m,M^m)$ is a solution of BSDE$^T(\xi,f_{m}+dV)$.
Letting  $m\rightarrow \infty$ in (\ref{eq2.4}) and repeating the
above argument, with obvious modification, shows the existence of
a solution of BSDE$^T(\xi,f+dV)$.
\end{dow}

\begin{stw}
\label{tw3.1234} Assume that  {{\rm(A1), (A2), (A4),
(A5)}} are satisfied and
\[
E\int_0^T|f(r,0)|\,dr<\infty.
\]
Then there exists a solution of {\rm BSDE}$^T(\xi,f+dV)$.
\end{stw}
\begin{dow}
Let $\xi_n=(\xi\wedge n)\vee (-n)$, $V^n_t=\int_{0}^{t\wedge
n}\mathbf{1}_{\{|V|_r\le n\}}\,dV_r$, and let
\[
f_n(t,y)=f(t,y)-f(t,0)+T_n(f(t,0))\cdot g_n(t),
\]
where $g_n(t)=1/(1+t^2/n)$. Observe that the data $\xi^n, V^n,
f_n$ satisfy the assumptions of Lemma \ref{lm3.123}. Therefore for
every $n\ge 1$ there exists a solution $(Y^n,M^n)$ of
BSDE$^T(\xi_n,f_n+dV^n)$. By the Tanaka-Meyer formula and (A2),
for $n<m$ we have
\begin{align*}
|Y^n_t-Y^m_t|&\le E\Big( |\xi_n-\xi_m|+\int_n^T\,d|V|_r
+\int_0^T\mathbf{1}_{\{n<|V|_r\le m\}}\,d|V|_r\\
&\quad+\int_0^T|T_n(f(r,0))g_n(r)-T_m(f(r,0))g_m(r)|\,dr\Big).
\end{align*}
By \cite[Lemma 6.1]{BDHPS},
\begin{equation}
\label{eq3.222}
E\sup_{t\ge 0}|Y^n_t-Y^m_t|^q\rightarrow 0.
\end{equation}
for every $q\in (0,1)$. It follows in particular that there is an
adapted c\`adl\`ag process  $Y$ such that $Y^n\rightarrow Y$ in
ucp. By the Tanaka-Meyer formula,
\[
|Y^{n}_t|\le E\big(|\xi|+\int_0^T|f(r,0)|\,dr +\int_0^T\,
d|V|_r|\FF_t\big)=:X_t.
\]
Furthermore,  by \cite[Lemma 2.3]{KR:JFA} and Fatou's lemma,
\begin{equation}
\label{eq3.333}
E\int_0^T|f(r,Y_r)|\,dr\le E\big(|\xi|
+\int_0^T|f(r,0)|\,dr+\int_0^T\, d|V|_r|\FF_t\big).
\end{equation}
Set
\[
\tau_k=\inf\{t\ge 0:\int_0^t|f(r,X_r)|\,dr\ge k\}.
\]
For every $a\ge 0$ we have
\[
Y^n_t=E\big(Y^n_{\tau_k\wedge a}+\int_t^{\tau_k\wedge a}
f_n(r,Y^{n}_r)\,dr+\int_t^{\tau_k\wedge a}\,dV^n_r|\FF_t\big),
\quad t\le\tau_k\wedge a.
\]
Letting $n\rightarrow\infty$ in the above equality and using (A4),
(A5) and (\ref{eq3.222}), (\ref{eq3.333}) we get
\[
Y_t=E\big(Y_{\tau_k\wedge a}+\int_t^{\tau_k\wedge a}f(r,Y_r)\,dr
+\int_t^{\tau_k\wedge a}\,dV_r|\FF_t\big),\quad t\le\tau_k\wedge a.
\]
Letting now $k,a\rightarrow \infty$ in the above equality and
using  (\ref{eq3.222}), (\ref{eq3.333}) we obtain
\[
Y_t=E\big(\xi+\int_t^{T}f(r,Y_r)\,dr+\int_t^{T}\,dV_r|\FF_t\big),
\quad t\le T.
\]
Set
\[
M_t=E\big(\xi+\int_0^Tf(r,Y_r)\,dr+\int_0^T\,dV_r|\FF_t\big)-Y_0,
\quad
t\le T.
\]
It is easily seen that the pair $(Y,M)$ is a solution of
BSDE$^T(\xi,f+dV)$.
\end{dow}

\begin{tw}
\label{stw3.0} Let  {{\rm(A1), (A2), (A4), (A5)}} be satisfied.
Assume also that there exists a semimartingale $S$ such that $S$
is a difference of supermartingales of class {\rm (D)} and
\[
E\int_0^T|f(r,S_r)|\,dr<\infty.
\]
Then there
exists a solution of {{\rm BSDE}}$^T(\xi,f+dV)$. Moreover,
\[
E\int_0^T|f(r,Y_r)|\,dr\le E\Big(|\xi|+|S_T|+\int_0^T|f(r,S_r)|\,dr +
E\int_0^T\,d|V|_r+\int_0^T\,d|C|_r\Big),
\]
where $S_t=S_0+C_t+N_t$ is the Doob-Meyer decomposition of $S$.
\end{tw}
\begin{dow}
Set
\[
f_S(t,y)=f(t,S_t+y),\quad \tilde{\xi}=\xi-S_T,\quad \tilde{V}_t=V_t-C_t.
\]
By Proposition \ref{tw3.1234} there exists a unique solution
$(\tilde{Y},\tilde{M})$ of {{\rm
BSDE}}$^T(\tilde{\xi},f_S+d\tilde{V})$. Set
$(Y,M)=(\tilde{Y}+S,\tilde{M}+N)$. Then $(Y,M)$ is a solution of
{{\rm BSDE}}$^T(\xi,f+dV)$. By \cite[Lemma 2.3]{KR:JFA},
\[
E\int_0^T|f_S(r,\tilde{Y}_r)|\,dr\le
E(|\tilde{\xi}|+\int_0^T|f_S(r,0)|\,dr+\int_0^T\,d|\tilde{V}|_r),
\]
which implies the desired inequality.
\end{dow}

\subsection{Systems of equations}

In the rest of this section we assume that $N\ge1$.

\begin{df}
We say that a pair $(Y,M)$ is a subsolution (resp. supersolution)
of BSDE$^T(\xi,f+dV)$ if there exist $\underline{\xi},
\underline{V}$ (resp. $\overline{\xi},\, \overline{V}$) such that
$\underline{\xi}\le \xi,\, d\underline{V}\le dV$ (resp.
$\overline{\xi}\ge \xi,\, d\overline{V}\ge dV$),
$E(\int_0^T\,d|\underline{V}|_r+|\underline{\xi}|)<\infty$ (resp.
$E(\int_0^T\,d|\overline{V}|_r+|\overline{\xi}|)<\infty$) and
$(Y,M)$ is a solution of
BSDE$^T(\underline{\xi},f+d\underline{V})$ (resp.
BSDE$^T(\overline{\xi},f+d\overline{V})$).
\end{df}

We will make the following assumption:
\begin{enumerate}
\item[(A6)] there exist a  subsolution
$(\underline{Y},\underline{M})$ and a
supersolution $(\overline{Y},\overline{M})$
of BSDE$^T(\xi,f+dV)$ such that
\[
\underline{Y}\le \overline{Y},\qquad
\sum_{j=1}^NE\Big(\int_0^T|f^j(r,\underline{Y}_r;S^j_r)|\,dr
+\int_0^T|f^j(r,\overline{Y}_r;S^j_r)|\,dr\Big)<\infty
\]
for some semimartingale $S$ which is a difference of supermartingales of class (D).
\end{enumerate}

\begin{prz}
\label{prz3.1} Let assumptions  (A1)--(A5) hold. If
$\overline{f},\underline{f}$ satisfy (A4), (A5) and
\begin{equation}
\label{eq2.2} \sum_{j=1}^NE(\int_0^T|\overline{f}^j(r,S^j_r)|\,dr
+\int_0^T|\underline{f}^j(r,S^j_r)|\,dr)<\infty,
\end{equation}
for some semimartingale $S$ which is a difference of supermartingales of class (D),
then (A6) is satisfied with $(\underline{Y}^j,\underline{M}^j)$,
$(\overline{Y}^j,\overline{M}^j)$ being solutions of
BSDE$^T(\xi^j,\underline{f}^j+dV^j)$ and
BSDE$^T(\xi^j,\overline{f}^j+dV^j)$, respectively.
\end{prz}

\begin{prz}
\label{prz2.123} Assume that (A1), (A4), (A5) are satisfied, $T$
is bounded and $f$ is Lipschitz continuous in $y$ uniformly in
$t$, i.e. there exists $L>0$ such that
\[
|f(t,y)-f(t,y')|\le L|y-y'|, \quad y,y'\in\BR^N.
\]
Then (A6) is satisfied by the pairs $(\overline{Y},
\overline{M})$, $(\underline{Y}, \underline{M})$ defined by
\[
\overline{Y}^1=\overline{Y}^2=\ldots=\overline{Y}^N,\qquad
\overline{M}^1=\overline{M}^2=\ldots=\overline{M}^N,
\]
\[
 \underline{Y}^1=\underline{Y}^2=\ldots=\underline{Y}^N,\qquad
\underline{M}^1=\underline{M}^2=\ldots=\underline{M}^N,
\]
where $(\overline{Y}^1,\overline{M}^1)$,
$(\underline{Y}^1,\underline{M}^1)$ are solutions of
BSDE$^T(\xi^1\vee\dots\vee\xi^N, f^1\vee\dots\vee
f^N+dV^1\vee\dots\vee dV^N)$ and
BSDE$^T(\xi^1\wedge\dots\wedge\xi^N, f^1\wedge\dots\wedge
f^N+dV^1\wedge\dots\wedge dV^N)$, respectively.
\end{prz}

\begin{tw}
\label{tw3.1} Let assumptions \mbox{\rm{(A1)--(A5)}} hold, and let
\mbox{\rm(A6)}  be satisfied with some processes $\underline{Y}\,,
\overline{Y}$. Then there exists a minimal solution $(Y,M)$ of
{\rm BSDE}$^{T}(\xi,f+dV)$ such that $\underline{Y}\le Y\le
\overline{Y}$. Moreover,
\begin{equation}
\label{eq3.1}
E\int_0^T|f(r,Y_r)|\,dr<\infty
\end{equation}
 and $M$ is a uniformly integrable martingale.
\end{tw}
\begin{dow}
Let $(\overline{Y},\overline{M})$, $(\underline{Y},\underline{M})$
be as in (A6). Let $Y^0:=\underline{Y}$ and $(Y^{n,j},M^{n,j})$,
$j=1,\dots,N$, be a solution of
BSDE$^T(\xi^j,f^j(\cdot,Y^{n-1};\cdot)+dV^j)$. Then
\begin{equation}
\label{eq3.2} Y_t^{n,j}=Y^{n,j}_{T\wedge a}+\int_t^{T\wedge a}
f^j(r,Y^{n-1}_r;Y^{n,j}_r)\,dr+\int_t^{T\wedge a} dV_r^j
-\int_t^{T\wedge a}\, dM^{n,j}_r,\quad t\in [0,T\wedge a].
\end{equation}
By \cite[Proposition 3.1]{KR:JFA},
\begin{equation}
\label{eq3.3} Y^n\le Y^{n+1},\qquad Y^n\le \overline{Y}.
\end{equation}
Therefore letting $n\rightarrow\infty$ in (\ref{eq3.2}) we get
\[
Y_t^{j}=Y^{j}_{T\wedge a}+\int_t^{T\wedge a} f^j(r,Y_r)\,dr
+\int_t^{T\wedge a} dV_r^j-\int_t^{T\wedge a}\, dM^{j}_r,
\quad t\in [0,T\wedge a],
\]
where $Y_t=\lim_{n\rightarrow\infty}Y^n_t$ and
$M_t=\lim_{n\rightarrow \infty}M^n_t,\, t\in [0,T\wedge a]$. The
process $M$ is a local martingale, because by (\ref{eq3.3}) the
sequence $\{M^n\}$ is locally uniformly integrable as all the
other terms in (\ref{eq3.2}) are locally uniformly integrable with
respect to $n$. To show that the pair $(Y,M)$ is a solution of
BSDE$^T(\xi,f+dV)$ it remains to prove  that $Y_{T\wedge
a}\rightarrow \xi$ as $a\rightarrow \infty$. If $T$ is finite,
this follows immediately from the fact that $Y^n_t\nearrow Y_t,\,
t\le T$. In general case an additional argument is required. By
Theorem \ref{stw3.0} there exists a solution
$(\overline{X}^j,\overline{N}^j)$ of
BSDE$^T(\xi,f^j(\overline{Y};\cdot)+dV^j)$ and a solution
$(\underline{X}^j,\underline{N}^j)$ of
BSDE$^T(\xi,f^j(\underline{Y};\cdot)+dV^j)$. Moreover, by
\cite[Proposition 3.1]{KR:JFA}, $\underline{X}_t\le Y^n_t\le
\overline{X}_t,\, t\in [0,T\wedge a],\, a\ge 0$, which  implies
the desired convergence. By (\ref{eq3.3}), (A6) and Theorem
\ref{stw3.0},
\begin{align*}
E\int_0^T|f^j(r,Y_r)|\,dr&\le E\Big(|\xi^j|+|S^j_T|
+\int_0^T|f^j(r,Y_r;S^j_r)|\,dr+\int_0^T\,d|V^j|_r
+\int_0^T\,d|C^j|_r\Big)\\& \le
E\Big(|\xi^j|+|S^j_T|+\int_0^T\,d|V^j|_r+\int_0^T\,d|C^j|_r\Big)\\
&\quad+ E\Big(\int_0^T|f^j(r,\underline{Y}_r;S_r^j)|\,dr
+\int_0^T|f^j(r,\overline{Y}_r;S^j_r)|\,dr\Big)<\infty.
\end{align*}
From this and the fact that $Y$ is of class (D) we conclude that
\[
M_t=E\big(\xi+\int_0^Tf(r,Y_r)\,dr+\int_0^T\,dV_r|\FF_t\big)
-Y_0,\quad t\in [0,T].
\]
It follows that $M$ is a uniformly integrable martingale.  Let
$(Y^*,M^*)$ be a solution of BSDE$^T(\xi,f+dV)$ such that
$\underline{Y}\le Y^*\le \overline{Y}$. Then by \cite[Proposition
3.1]{KR:JFA}, $Y^n\le Y^*$, $n\ge 0$, which implies that $Y\le
Y^*$.
\end{dow}

\begin{wn}
\label{wn3.1} Assume that the data $(\xi,f,V)$, $(\xi',f',V')$
satisfy  \mbox{\rm(A1)--(A5)}, and that \mbox{\rm(A6)} is
satisfied with the same processes $\underline{Y}\,,\overline{Y}$.
Moreover, assume that
\[
\xi\le\xi',\quad f\le f',\quad dV\le dV',
\]
and that  $(Y,M)$ (resp. $(Y',M')$) is the minimal solution of
BSDE$^T(\xi,f+dV)$ (resp. BSDE$^T(\xi',f'+dV')$) such that
$\underline{Y}\le Y\le \overline{Y}$ (resp. $\underline{Y}\le
Y'\le \overline{Y}$). Then
\[
Y_t\le Y'_t\quad t\in [0,T].
\]
\end{wn}
\begin{dow}
Follows from the construction of processes $Y,Y'$ (see Theorem
\ref{tw3.1}) and \cite[Proposition 3.1]{KR:JFA}.
\end{dow}

\section{Systems of BSDEs with oblique reflection} \label{sec3}

Consider a family $\{h_{j,i}; i,j=1,\dots,N\}$ of measurable
functions $h_{j,i}:\Omega\times\BR^+\times\BR\rightarrow \BR$ such
that $h_{j,i}(\cdot,y^i)$ is progressively measurable for every
$y^i\in\BR$. For given sets  $A_j\subset
\{1,\dots,j-1,j+1,\dots,N\}$, $j=1,\dots,N$, set
\[
H^j(t,y)=\max_{i\in A_j}h_{j,i}(t,y^i), \quad
H(t,y)=(H^1(t,y),\dots,H^N(t,y)),\quad t\in\BR^+,\, y\in\BR^N.
\]
We adopt the convention that the maximum over the empty set equals
$-\infty$. Consequently, if $A_j=\emptyset$ for some $j$, then
$H^j(t,y)=-\infty$.

Apart from (A1)--(A6) we will also need the following assumptions:
\begin{enumerate}
\item[(A7)] There exist a subsolution $(\underline{Y},\underline{M})$ and
 a supersolution $(\overline{Y},\overline{M})$
of BSDE$^T(\xi,f+dV)$  such that
\[
H(\cdot,\overline{Y})\le \overline{Y},\quad \underline{Y}\le
\overline{Y},\quad
\sum_{j=1}^NE\big(\int_0^T|f^j(r,\underline{Y}_r;\overline{Y}^j_r)|\,dr
+\int_0^T|f^j(r,\overline{Y}_r)|\,dr\Big)<\infty,
\]
\item[(A8)] $(t,y)\mapsto H^j(t,y)$
is continuous, $y\mapsto H^j(t,y)$ is nondecreasing and
\[
\limsup_{(t,y)\rightarrow (\infty,\xi)} H^j(t,y)\le\xi^j.
\]
\end{enumerate}

\begin{prz}
Let assumptions  (A1)--(A5) hold. Moreover, assume that
$\overline{f},\underline{f}$ satisfy (A4), (A5), (\ref{eq2.2}) (with $S^1=\ldots=S^N$),
and  $h_{j,i}(t,a)\le a$ for every $a\in\BR$. Let
\[
\overline{Y}^1=\overline{Y}^2=\ldots=\overline{Y}^N,\qquad
\overline{M}^1=\overline{M}^2=\ldots=\overline{M}^N,
\]
where $(\overline{Y}^1,\overline{M}^1)$ is a  solution of
BSDE$^T(\sum_{j=1}^N\xi^{j,+},
\sum_{j=1}^N(\overline{f}^{j,+}+dV^{j,+}))$. By
$(\underline{Y}^j,\underline{M}^j)$ denote a solution of
BSDE$^T(\xi^j,\underline{f}^j+dV^j)$. The solutions
$(\overline{Y}^1,\overline{M}^1)$,
$(\underline{Y}^j,\underline{M}^j)$ exist by Theorem \ref{stw3.0}.
By \cite[Proposition 3.1]{KR:JFA}, $\underline{Y}\le\overline{Y}$.
It follows that the pair $(\underline{Y},\overline{Y})$  satisfies
(A7).
\end{prz}

\begin{prz}
Let assumptions of Example \ref{prz2.123} hold, and let
$h_{j,i}(t,a)\le a$ for every $a\in\BR$. Then the processes
$(\overline{Y},\overline{M})$, $(\underline{Y},\underline{M})$
defined in Example \ref{prz2.123} satisfy (A7).
\end{prz}

\begin{uw}
\label{uw.121} By Theorem \ref{stw3.0}, if $f$ satisfies (A2)
and (A7) then
\begin{align*}
E\int_0^T|f(r,\underline{Y}_r)|\, dr &\le
E\Big(|\underline{\xi}|+|\overline{\xi}|+\sum_{j=1}^N\int_0^T|f^j(r,\underline{Y}_r;\overline{Y}^j_r)|\,dr\\&\quad +
\sum_{j=1}^N\int_0^T|f^j(r,\overline{Y}_r)|\,dr
+\int_0^T\, d|\underline{V}|_r+\int_0^T\, d|\overline{V}|_r\Big).
\end{align*}
\end{uw}

\subsection{Existence of solutions}

\begin{df}
We say that a triple $(Y,M,K)$ of adapted c\`adl\`ag  processes is
a solution of BSDE with oblique reflection (\ref{eq1.3}) if  $Y$
is of class (D), $M$ is a local martingale with $M_0=0$, $K$ is an
increasing process with $K_0=0$ and (\ref{eq1.3}) is satisfied.
\end{df}

If $A_j=\emptyset$, then by convention, $H^j=-\infty$, and hence
$Y^j$ has no lower barrier. We then take $K^j=0$ in the above
definition.

Recall the following definition from \cite{K:SPA}.
\begin{df}
\label{def2.11} Let $N=1$, and $L$ let be a c\`adl\`ag process. We
say that a triple $(Y,M,K)$ of adapted c\`adl\`ag processes is a
solution of reflected BSDE on the interval $[0,T]$ with terminal
condition $\xi$, right-hand side $f+dV$ and lower barrier $L$
(RBSDE$^T(\xi,f+dV,L)$ for short) if
\begin{enumerate}
\item[(i)] $Y$ is of class (D), $M$ is a local martingale
with $M_0=0$, $K$ is an increasing process with $K_0=0$,
\item[(ii)] $Y_t\ge L_t,\, t\in [0,T\wedge a]$,
$\int_0^{T\wedge a}(Y_{t-}-L_{t-})\,dK_t=0$ for every $a\ge 0$,
\item[(iii)] $\int_0^{T\wedge a}|f(t,Y_t)|\,dt<\infty,\quad a\ge 0,$
\item[(iv)] for every $a\ge 0$,
\[
Y_t=Y_{T\wedge a}
+\int_t^{T\wedge a}f(r,Y_r)\,dr+\int_t^{T\wedge a}\,dV_r
+\int_t^{T\wedge a}\,dK_r-\int_t^{T\wedge a}\,dM_r,
\quad t\in [0,T\wedge a],
\]
\item[(v)] $Y_{T\wedge a}\rightarrow \xi$ $P$-a.s. as $a\rightarrow\infty$.
\end{enumerate}
\end{df}

Observe that  a triple $(Y,M,K)$  is a solution of (\ref{eq1.3})
if and only if $(Y^j,M^j,K^j)$ is a solution of
RBSDE$^T(\xi^j,f^j(\cdot,Y;\cdot)+dV^j,H^j(\cdot,Y))$ for every
$j=1,\dots,N$.

\begin{uw}
\label{uw3.rev1} If (\ref{eq2.01}) is satisfied then
$EK_T<\infty$, $M$ is a uniformly integrable martingale and
\[
Y_t=\xi +\int_t^{T}f(r,Y_r)\,dr+\int_t^{T}\,dV_r
+\int_t^{T}\,dK_r-\int_t^{T}\,dM_r,\quad t\in [0,T].
\]
Indeed, localizing the local martingale $M$ we easily deduce that
$EK_T<\infty$. The remaining two assertions then follow from
Remark \ref{uw3.0}.
\end{uw}

\begin{uw}
\label{uw3.rev2} Let $(Y,M,K)$ be a solution of
RBSDE$^T(\xi,f+dV,L)$. Under the assumptions of Remark
\ref{uw3.rev1},
\begin{equation}
\label{eq.rev3.1} Y_t=\esssup_{\tau\ge t}E\Big(\int_t^{T\wedge
\tau} f(r,Y_r)\,dr+\int_t^{T\wedge \tau} \,dV_r+
L_\tau\mathbf{1}_{\{ \tau<T\}}
+\xi\mathbf{1}_{\{T\wedge\tau=T\}}|\FF_t\Big).
\end{equation}
To see this, we first observe that by Remark \ref{uw3.rev1}, for
every stopping time $\tau\ge t$,
\begin{align*}
Y_t&=E\Big(Y_{T\wedge\tau }+\int_t^{T\wedge \tau}
f(r,Y_r)\,dr+\int_t^{T\wedge\tau}\,dK_r+\int_t^{T\wedge \tau}
\,dV_r|\FF_t\big)\\
& \ge E\Big(\int_t^{T\wedge \tau} f(r,Y_r)\,dr+\int_t^{T\wedge
\tau} \,dV_r+ L_\tau\mathbf{1}_{\{\tau<T\}}
+\xi\mathbf{1}_{\{T\wedge\tau=T\}}|\FF_t\Big).
\end{align*}
This shows that $Y_t$ is greater then or equal to the right-hand
side of  (\ref{eq.rev3.1}). To get the opposite inequality, we
consider the stopping time
\[
D^\varepsilon_t=\inf\{s\ge t,\, L_s+\varepsilon\ge Y_s\}\wedge T.
\]
By the minimality property of $K$,
\begin{align*}
Y_t&=E(Y_{D^\varepsilon_t}+\int_t^{D^\varepsilon_t}
f(r,Y_r)\,dr+\int_t^{D^\varepsilon_t} \,dV_r|\FF_t)\\& \le
E(L_{D^\varepsilon_t}\mathbf{1}_{\{D^\varepsilon_t<T\}}
+\xi\mathbf{1}_{\{D^\varepsilon_t=T\}} +\int_t^{D^\varepsilon_t}
f(r,Y_r)\,dr+\int_t^{D^\varepsilon_t} \,dV_r|\FF_t)+\varepsilon,
\end{align*}
from which it follows that $Y_t$ is less then or equal to the
right-hand side of (\ref{eq.rev3.1}).
\end{uw}

In \cite{K:SPA} an existence result for RBSDE$^T(\xi,f+dV,L)$ is
proved under the assumption that $T$ is bounded. The general case
requires some modification of the proof given in \cite{K:SPA}. In
this modified proof we will need the following lemma.

\begin{lm}
\label{lm.rev1} Assume that $L^+$ is of class {\rm(D)},
$E|\xi|<\infty$  and $\limsup_{a\rightarrow \infty}L_{T\wedge a}
\le \xi$. Then
\begin{equation}
\label{eq3.04} \limsup_{a\rightarrow \infty}Y_{T\wedge a}\le \xi,
\end{equation}
where
\begin{equation}
\label{eq3.03} Y_t=\esssup_{\tau\ge t}
E(L_{\tau}\mathbf{1}_{\{\tau<T\}}+\xi\mathbf{1}_{\{T\wedge\tau=T\}}|\FF_t).
\end{equation}
\end{lm}
\begin{dow}
From the definition of $Y$ it follows that $Y_{t}=Y_{T\wedge t}$.
Therefore the assertion of the lemma is clear if $T<\infty$. 
Let $\varepsilon>0$. By the assumptions of the lemma, for a.e.
$\omega\in\Omega$ there exists $t_\omega$ such that
\[
L_t(\omega)\le\xi(\omega)+\varepsilon,\quad t\ge t_\omega.
\]
Let
\[
\Lambda_n=\{\omega\in\Omega;\, t_\omega\ge n\}.
\]
It is clear that $\Lambda_{n+1}\subset \Lambda_n$ and
$P(\bigcap_{n\ge 1} \Lambda_n)=0$. Since $L^+$ is of class (D),
 there is $\delta>0$ such that
if $A\in\FF$ and $P(A)<\delta$ then
$\sup_{\tau}\int_{A}(L^+_\tau\mathbf{1}_{\{T\wedge\tau<T\}}+|\xi|)\le
\varepsilon$. Choose  $N\in\mathbb{N}$ so that  $P(\Lambda_N)\le
\delta$. Then for $t\ge N$,
\begin{align*}
Y_t&=\esssup_{\tau\ge t} E((L_{\tau}\mathbf{1}_{\{\tau<T\}}
+\xi\mathbf{1}_{\{T\wedge\tau=T\}})
\mathbf{1}_{\Lambda_N^c}|\FF_t)\\& \quad+\esssup_{\tau\ge
t}E((L_{\tau} \mathbf{1}_{\{\tau<T\}}
+\xi\mathbf{1}_{\{T\wedge\tau=T\}})\mathbf{1}_{\Lambda_N}|\FF_t)\\&
\le \varepsilon +E(\xi\mathbf{1}_{\Lambda_N^c}|\FF_t)
+\esssup_{\tau\ge t}E((L^+_{\tau}\mathbf{1}_{\{\tau<T\}}
+|\xi|\mathbf{1}_{\{T\wedge\tau=T\}})\mathbf{1}_{\Lambda_N}|\FF_t)\\&
\le 2 \varepsilon +E(\xi\mathbf{1}_{\Lambda_N^c}|\FF_t).
\end{align*}
Letting $t\rightarrow \infty$ and then $N\rightarrow \infty$ we
get
$\limsup_{t\rightarrow\infty}Y_t
\le2\varepsilon +\xi$, which implies (\ref{eq3.04}).
\end{dow}

\begin{tw}
\label{tw.rev1} Let $N=1$. Assume that {\rm (A1), (A2), (A4),
(A5)} are satisfied and $L$ is a c\`adl\`ag adapted process such
that $\limsup_{a\rightarrow \infty}L_{T\wedge a}\le\xi$ and $L\le
X$ for some semimartingale $X$ such that $X$ is a difference of
supermartingales of class {\rm (D)} and
\[
E\int_0^T|f(r,X_r)|\,dr<\infty.
\]
Then there exists a solution $(Y,M,K)$ of RBSDE$^T(\xi,f+dV,L)$.
Moreover,
\[
E\int_0^T|f(r,Y_r)|\,dr+EK_T<\infty
\]
and $M$ is a uniformly integrable martingale.
\end{tw}
\begin{dow}
The proof runs as the proof of \cite[Theorem 2.13]{K:SPA}, with
small modifications. By Theorem \ref{stw3.0} there exists a
solution $(Y^n,M^n)$ of BSDE$^T(\xi,f_n+dV)$ with
\[
f_n(t,y)=f(t,y)+n(y-L_t)^-.
\]
By \cite[Proposition 3.1]{KR:JFA}, $Y^n\le Y^{n+1}$. As in
\cite{K:SPA} we construct a supersolution
$(\overline{X},\overline{N})$ of BSDE$^T(\xi,f+dV)$ such that
$\overline{X}\ge L$ and
\begin{equation}
\label{eq.rev1}
Y^1\le Y^n\le \overline{X},\quad n\ge 1.
\end{equation}
By Theorem \ref{stw3.0},
\begin{equation}
\label{eq.rev2}
E\int_0^T|f(r,Y^1_r)|\,dr+E\int_0^T|f(r,\overline{X}_r)|\,dr<\infty.
\end{equation}
Therefore by (A2), (\ref{eq.rev1}) and the Lebesgue dominated
convergence theorem,
\begin{equation}
\label{eq.rev3}
E\int_0^T|f(r,Y^n_r)-f(r,Y_r)|\,dr\rightarrow \infty,
\end{equation}
where $Y_t=\sup_{n\ge 1} Y^n_t,\, t\ge 0$. Repeating now, on each
interval $[0,T\wedge a]$, the reasoning following (2.22) in the
proof of \cite[Theorem 2.13]{K:SPA} we show that $Y$ is c\`adl\`ag
and there exists a predictable c\`adl\`ag increasing process $K$
with $K_0=0$ and a local martingale $M$ with $M_0=0$ such that for
every $a\ge0$,
\[
Y_t=Y_{T\wedge a}+\int_t^{T\wedge a} f(r,Y_r)\,dr
+\int_t^{T\wedge a}\,dV_r+\int_t^{T\wedge a}\,dK_r
-\int_t^{T\wedge a}\,dM_r,\quad t\in [0,T\wedge a]
\]
and
\[
Y\ge L,\quad \int_0^{T\wedge a}(Y_{r-}-L_{r-})\,dK_r=0.
\]
By (\ref{eq.rev1}), $Y$ is of class (D), which combined with
(\ref{eq.rev2}) yields $E\int_0^T|f(r,Y_r)|\,dr+EK_T<\infty$. This
inequality implies that $M$ is a uniformly integrable martingale
(see Remark \ref{uw3.rev1}). What is left is to show that
$Y_{T\wedge a}\rightarrow \xi,\, a\rightarrow \infty$.  By
(\ref{eq.rev1}) $\xi\le \liminf_{a\rightarrow \infty}Y_{T\wedge
a}$, so it suffices to show that
\begin{equation}
\label{eq3.05} \limsup_{a\rightarrow \infty}Y_{T\wedge a}\le \xi.
\end{equation}
Observe that the triple $(Y^n,M^n,K^n)$, where $
K^n_t=\int_0^t(Y^n_r-L_r)^-\,dr$, is a solution of
RBSDE$^T(\xi,f+dV,L^n)$ with $L^n_t=L_t-(Y^n_t-L_t)^{-}$.
Therefore by Remark \ref{uw3.rev2} and the definition of $L^n$,
\[
Y^n_t\le \esssup_{\tau\ge t}E\Big(\int_t^{T\wedge \tau}
f(r,Y^n_r)\,dr+\int_t^{T\wedge \tau} \,dV_r+
L_\tau\mathbf{1}_{\{T\wedge \tau<T\}}
+\xi\mathbf{1}_{\{T\wedge\tau=T\}}|\FF_t\Big).
\]
Letting $n\rightarrow$ and using  (\ref{eq.rev3}) we get
\[
Y_t\le \esssup_{\tau\ge t}E\Big(\int_t^{T\wedge \tau}
f(r,Y_r)\,dr+\int_t^{T\wedge \tau} \,dV_r+
L_\tau\mathbf{1}_{\{T\wedge \tau<T\}}
+\xi\mathbf{1}_{\{T\wedge\tau=T\}}|\FF_t\Big).
\]
From this and  Lemma \ref{lm.rev1} we conclude that (\ref{eq3.05})
is satisfied.
\end{dow}

To prove the existence result for (\ref{eq1.3}) we will need the
monotone convergence theorem for BSDEs stated below. In the case
of Brownian filtration this result was proved
in \cite{K:BSM,PengXu}. In the case of  general filtration it
follows from \cite{K:SPA}.

\begin{stw}
\label{stw3.1}
Let $N=1$ and  {{\rm (A1)}} be
satisfied. Assume that $(Y^n,M^n)$ is a solution of
{\rm BSDE}$^T(\xi,dV^n+dK^n)$, where $K^n$ is an increasing predictable
c\`adl\`ag process such that $K^n_0=0$, and $V^n$ is a finite
variation c\`adl\`ag process with $V^n_0=0$. Moreover, assume that
$Y^n\le Y^{n+1}$, there exists a c\`adl\`ag process $\overline{Y}$
of class {{\rm D}} such that $Y^n\le \overline{Y}$, and that
$\{|V^n|\}$ is locally bounded in $L^2$ and $V^n\rightarrow V$ in
ucp for some finite variation c\`adl\`ag process $V$. Then there
exists a local martingale $M$ with $M_0=0$ and a predictable
c\`adl\`ag increasing process $K$ with $K_0=0$ such that
for every $a\ge 0$,
\begin{align*}
Y_{t}=Y_{T\wedge a}+\int_t^{T\wedge a}\,dV_r+\int_t^{T\wedge a} dK_r-\int_{t}^{T\wedge
a}\,dM_r,\quad t\in [0,T\wedge a],
\end{align*}
where $Y_t=\sup_{n\ge 1}Y^n_t,\, t\in [0,T\wedge a],\,
a\ge 0$.
Moreover, if $T<\infty$ then the pair
$(Y,M)$ is a solution of {\rm BSDE}$^T(\xi,dV+dK)$.
\end{stw}
\begin{dow}
It is enough to repeat the arguments between (2.22)--(2.28) in the
proof of \cite[Theorem 2.13]{K:SPA} with
$\overline{X}=\overline{Y}$ and with $(\int_0^\cdot
f(r,Y^n_r)\,dr, V^n)$ replaced by $(V^n,0)$.
\end{dow}

\begin{tw}
\label{tw3.2} Let assumptions {{\rm (A1)--(A5), (A8)}} hold, and
let \mbox{\rm(A7)}  be satisfied with some processes
$\underline{Y}\,,\overline{Y}$. Then there exists a minimal
solution $(Y,M,K)$ of \mbox{\rm(\ref{eq1.3})}  such that
$\underline{Y}\le Y\le \overline{Y}$.
\end{tw}
\begin{dow}
Let $(Y^0,M^0,K^0):=(\underline{Y},\underline{M},0)$. We define
$(Y^{j,n},M^{j,n},K^{j,n})$ to be a solution of RBSDE$^T(\xi^j,
f^j(\cdot,Y^{n-1};\cdot)+dV,H^j(\cdot,Y^{n-1}))$ (see Definition
\ref{def2.11}). It exists by Theorem \ref{tw.rev1}. For every
$a\ge 0$ we have
\begin{equation}
\label{eq3.4} \left\{
\begin{array}{l}
Y^{n,j}_{t}=Y^{n,j}_{T\wedge a}+\int_{t}^{T\wedge a}
f^{j}(r,Y^{n-1}_r;Y^{n,j}_r)\,dr+\int_t^{T\wedge a}\,dV^j_r\\
\qquad\qquad\qquad\qquad\qquad\qquad+\int_t^{T\wedge a}
dK^{n,j}_r-\int_{t}^{T\wedge a}\,dM_r^{n,j},
\medskip\\
Y^{n,j}_{t}\ge H^j(t,Y^{n-1}_t),\quad t\in [0,T\wedge a],
\medskip\\
\int_0^{T\wedge a}
(Y^{n,j}_{t-}-H^j_{t-}(\cdot,Y^{n-1}))\,dK^{n,j}_{t}=0.
\end{array}
\right.
\end{equation}
Moreover by (A2), (A3), (A8) and \cite[Proposition 2.1]{K:SPA},
\begin{equation}
\label{eq3.5}
Y^n\le Y^{n+1}\le\overline{Y},\quad n\ge 0.
\end{equation}
By Proposition \ref{stw3.1} there exists an increasing predictable
c\`adl\`ag process $K$ with $K_0=0$ and a local martingale $M$
with $M_0=0$ such that for every $a\ge 0$,
\begin{align}
\label{eq3.6} Y^{j}_{t}&=Y^{j}_{T\wedge a}+\int_{t}^{T\wedge a}
f^{j}(r,Y_r)\,dr+\int_t^{T\wedge a}\,dV^j_r \nonumber\\
&\quad+\int_t^{T\wedge a} dK^{j}_r-\int_{t}^{T\wedge
a}\,dM_r^{j},\quad t\in [0,T\wedge a],
\end{align}
where $Y_t=\sup_{n\ge 0} Y^n_t$. By (\ref{eq3.4}) and (A8) we also
have $Y^j\ge H^j(\cdot,Y)$. Let
$(\overline{X}^j,\overline{N}^j,\overline{K}^j)$ denote a solution
of RBSDE$^T(\xi^j,f^j(\overline{Y};\cdot)+dV^j,
H^j(\cdot,\overline{Y}))$ and
$(\underline{X}^j,\underline{N}^j,\underline{K}^j)$ denote a
solution of
RBSDE$^T(\xi^j,f^j(\underline{Y};\cdot)+dV^j,H^j(\cdot,\underline{Y}))$.
By (\ref{eq3.5}) and \cite[Proposition 2.1]{K:SPA},
$\underline{X}_t\le Y^n_t\le \overline{X}_t,\, t\in [0,T\wedge
a],\, a\ge 0$. This implies that $Y_{T\wedge a}\rightarrow \xi$ as
$a\rightarrow \infty$. What is left is to show that $K$ satisfies
the minimality condition. Set
\[
\tau_k=\inf\Big\{t\ge 0:
\sum_{j=1}^N\int_0^t|f^j(r,\overline{Y}_r;\underline{Y}^j_r)|
+|f^j(r,\underline{Y}_r;\overline{Y}^j_r)|\,dr\ge k\Big\}\wedge T.
\]
Then on $[0,\tau_k]$
\begin{equation}
\label{eq2.09} Y^{j}_t=\esssup_{t\le\tau}
E\Big(\int_t^{\tau_k\wedge\tau} f^j(r,Y_r)\,dr+\int_t^{\tau_k\wedge\tau}\,dV_r
+H^j(\tau,Y_{\tau})\mathbf{1}_{\{\tau<\tau_k\}}
+Y^{j}_{\tau_k}\mathbf{1}_{\{\tau_k\wedge\tau=\tau_k\}}|\FF_t\Big).
\end{equation}
Indeed, by Remark \ref{uw3.rev2},
\begin{align*}
Y^{n,j}_t&=\esssup_{t\le\tau}E\Big(\int_t^{\tau_k\wedge\tau}
f^j(r,Y^{n-1}_r;Y^{n,j}_r)\,dr+\int_t^{\tau_k\wedge\tau}\,dV_r\\
&\qquad\qquad\qquad\qquad
+H^j(\tau,Y^{n-1}_{\tau})\mathbf{1}_{\{\tau<\tau_k\}}
+Y^{n,j}_{\tau_k}\mathbf{1}_{\{\tau_k\wedge\tau=\tau_k\}}|\FF_t\Big),
\end{align*}
so by (A3) and (A8),
\begin{align*}
Y^{n,j}_t&\le\esssup_{t\le\tau}E\Big(\int_t^{\tau_k\wedge\tau}
f^j(r,Y_r;Y^{n,j}_r)\,dr+\int_t^{\tau_k\wedge\tau}\,dV_r\\
&\qquad\qquad\qquad\qquad
+H^j(\tau,Y_{\tau})\mathbf{1}_{\{\tau<\tau_k\}}
+Y^{j}_{\tau_k}\mathbf{1}_{\{\tau_k\wedge\tau=\tau_k\}}|\FF_t\Big).
\end{align*}
Letting $n\rightarrow \infty$ and using (A4) we see that $Y^j$ is
less then or equal to the right-hand side of (\ref{eq2.09}). The
opposite inequality follows from the fact that the process
$Y^j+\int_0^\cdot f(r,Y_r)\,dr+\int_0^\cdot\,dV_r$ is a
supermartingale  which dominates the process $L=\int_0^\cdot
f(r,Y_r)\,dr +\int_0^\cdot\,dV_r
+H^j(\cdot,Y)\mathbf{1}_{\{\cdot<\tau_k\}}
+Y^j_{\tau_k}\mathbf{1}_{\{\cdot=\tau_k\}}$. Thus (\ref{eq2.09})
is proved. By (\ref{eq2.09}) and Remark \ref{uw3.rev2},
\[
\int_0^{\tau_k} (Y^{j}_{t-}-H^j_{t-}(\cdot,Y))\,dK^{j}_{t}=0.
\]
Letting $k\rightarrow\infty$  gives the above inequality on every
interval $[0,T\wedge a],\, a\ge 0$. Let $(Y^*,M^*,K^*)$ be a
solution of (\ref{eq1.3}) such that $\underline{Y}\le Y^*\le
\underline{Y}$. By \cite[Proposition 2.1]{K:SPA}, $Y^n\le Y^*$,
$n\ge 0$. Hence $Y\le Y^*$.
\end{dow}

\begin{uw}
\label{uw3.1} If $K^n, K, V$ from the proof of Theorem \ref{tw3.2}
are continuous, then $Y^n\nearrow Y,\, K^n\rightarrow K$ in ucp.
Indeed, in this case
\begin{equation}
\label{eq2.7} ^pY^n_t=Y^n_{t-}\,,\qquad ^pY_t=Y_{t-},
\end{equation}
where $^pY^n$, $^pY$ denote predictable projections of $Y^n$ and
$Y$, respectively. It is known that $Y^n\nearrow Y$ implies that
${}^pY^n\nearrow {}^pY$. By this and (\ref{eq2.7}),
$Y^n_{t-}\nearrow Y_{t-}$, $t\in [0,T\wedge a]$, $a\ge 0$.
Therefore by the generalized Dini theorem (see \cite[p.
185]{DM2}), $Y^n\nearrow Y$ in ucp. The convergence of $\{K^n\}$
now follows from \cite[Theorem 1.8]{Jacod} (for details see the
reasoning at the beginning of page 4220 in \cite{K:SPA}).
\end{uw}

\begin{uw}
\label{uw3.2} If in Theorem \ref{tw3.2} we assume additionally
that
\begin{equation}
\label{eq3.7}
\sum_{i=1}^NE\Big(\int_0^T|f^j(r,\overline{Y}_r;\underline{Y}^j_r)|\,dr\Big)<\infty,
\end{equation}
where $\overline{Y}$, $\underline{Y}$ are processes from  (A7),
then $E\int_0^T|f(r,Y_r)|\,dr+E\int_0^T\,dK_r<\infty$ and $M$ is a
uniformly integrable martingale. This follows immediately from
(A2), (A3) and the fact that $\underline{Y}\le Y\le \overline{Y}$.
\end{uw}

\begin{uw}
\label{uw3.3} In Theorem \ref{tw3.2} assume additionally that $h$
is strictly increasing with respect to $y$ ((A8) implies only that
it is nondecreasing), and the following condition considered in
\cite{HZ}:
\begin{enumerate}
\item[(A9)] there are no $(y_1,\dots, y_k)\in\BR^k$ and
$j_2\in A_{j_1},\dots,j_k \in A_{j_{k-1}}, j_1\in A_{j_k}$ such
that
\[
y_1=h_{j_1,j_2}(t,y_2),\,\, y_2=h_{j_2,j_3}(t,y_3),\dots, y_{k-1}
=h_{j_{k-1},j_k}(t,y_k),\,\, y_k=h_{j_k,j_1}(t,y_1)
\]
\end{enumerate}
Moreover, assume that the underlying filtration is quasi-left
continuous and $V$ is continuous. Then $K$ is continuous. Indeed,
since the filtration is quasi-left continuous and $V$ is
continuous, $\Delta K_{\tau}=-\Delta Y_\tau$ for every predictable
stopping time $\tau$. Therefore in the same way as in Step 4 of
the proof of \cite[Theorem 3.2]{HZ} one can  show that $\Delta
K_\tau=0$. Since $K$ is predictable and $\tau$ is an arbitrary
predictable stopping time, applying the predictable cross-section
theorem (see \cite[Theorem 86, p. 138]{DM1}) shows that $K$ is
continuous.
\end{uw}

\subsection{Approximation via penalization}

Let us consider the following system of BSDEs
\begin{align}
\label{eq3.8} Y^{j,n}_t&=\xi^j+\int_t^T
f^j(r,Y^n_r)\,dr+\int_t^T\,dV^j_r\nonumber\\
&\quad+\int_t^Tn(Y^{j,n}_r-H^j(r,Y^n_r))^-\,dr-\int_t^T\,dM^{j,n}_r.
\end{align}

\begin{tw}
\label{tw3.3} Let \mbox{\rm (A1)--(A5), (A8)}   hold, and let  \mbox{\rm(A7)} be satisfied  with
some processes $\underline{Y}\,,\overline{Y}$. Then there exists a
minimal solution $(Y^n,M^n)$ of {{\rm (\ref{eq3.8})}} such that
$\underline{Y}\le Y^n\le\overline{Y}$. Moreover, $Y^n_t\nearrow
Y_t,\, t\in [0,T\wedge a],\, a\ge 0$, where $(Y,M,K)$ is the
minimal solution of {{\rm (\ref{eq1.3})}} such that
$\underline{Y}\le Y\le \overline{Y}$.
\end{tw}
\begin{dow}
Observe that by (A7), $(\overline{Y},
\overline{M})$ is a supersolution of (\ref{eq3.8}) such that
\[
E\int_0^T|f(r,\overline{Y}_r)|\,dr<\infty.
\]
It is clear that $(\underline{Y},\underline{M})$ is a subsolution
of (\ref{eq3.8}),  and that, by (A7),
\[
\sum_{i=1}^{N}E\int_0^T|f^j(r,\underline{Y}_r;\overline{Y}^j_r)|\,dr<\infty.
\]
Hence (A6) is satisfied for equation (\ref{eq3.8}) with
$\underline{Y},\overline{Y}$ and with $S=\overline{Y}$. Since the
other assumptions of Theorem \ref{tw3.1} for equation
(\ref{eq3.8}) are also  satisfied, there exists a minimal solution
$(Y^n,M^n) $ of (\ref{eq3.8}) such that $\underline{Y}\le Y^n\le
\overline{Y}$. By Corollary \ref{wn3.1}, $Y^n\le Y^{n+1}$.
Therefore repeating step by step the arguments from the proof of
Theorem \ref{tw3.2} (see also the end of the proof of Theorem \ref{tw.rev1})
we show that there exists a local martingale
$\tilde{M}$ and
an  increasing c\`adl\`ag process $\tilde{K}$ such that the triple
$(\tilde{Y},\tilde{M},\tilde{K})$, where
$\tilde{Y}_t=\lim_{n\rightarrow \infty}Y^n_t$, $t\in [0,T\wedge
a]$, $a\ge 0$, is a solution of (\ref{eq1.3}). What is now left is
to show  that $\tilde{Y}=Y$, where $(Y,M,K)$ is the minimal
solution of (\ref{eq1.3}) such that $\underline{Y}\le Y\le
\overline{Y}$. Of course, $Y\le \tilde{Y}$. Moreover, since
$Y^j\ge H^j(\cdot,Y)$, we have
\begin{align*}
Y^{j}_t&=\xi^j+\int_t^T
f^j(r,Y_r)\,dr+\int_t^T\,dV^j_r+\int_t^T\,dK^j_r\\&\quad
+\int_t^Tn(Y^{j}_r-H^j(r,Y_r))^-\,dr-\int_t^T\,dM^{j}_r.
\end{align*}
By this and  Corollary \ref{wn3.1}, $Y^n\le Y$. Hence
$\tilde{Y}\le Y$, which completes the proof.
\end{dow}

\begin{uw}
Set
\[
K^{n,j}_t=\int_0^tn(Y^{n,j}_r-H^j(r,Y^n_r))^-\,dr.
\]
If the processes $K, V$ from  Theorem \ref{tw3.2} are continuous,
then $Y^n\nearrow Y$ and $K^n\rightarrow K$ in ucp. This follows
by the same method as in Remark \ref{uw3.1}.
\end{uw}

\section{Switching problem}
\label{sec4}

In what follows by a strategy we mean a pair $\mathcal{S}=
(\{\xi_n,\, n\ge 1\},\, \{\tau_n,\, n\ge 1\})$, where $\{\tau_n,\,
n\ge 1\}$ is an increasing  sequence of stopping times such that
\[
P(\tau_n<T,\,\forall\,\, n\ge 1)=0,
\]
and $\{\xi_n,\,n\ge 1\}$ is a sequence of random variables taking
values in $\{1,\dots, N\}$ such that $\xi_n$ is
$\FF_{\tau_n}$-measurable for each $n\in\BN$.

The set of all strategies we denote by $\mathbf{A}$. By
$\mathbf{A}_t$ we denote the set of all strategies $\mathcal{S}\in
\mathbf{A}$ such that $\tau_1\ge t$. For
$\mathcal{S}\in\mathbf{A}$ we set
\[
w^j_t=j\mathbf{1}_{[0,\tau_1)}(t) +\sum_{n\ge 1}\xi_n
\mathbf{1}_{[\tau_{n},\tau_{n+1})}(t).
\]

\begin{uw}
\label{uw3.4} Let $L$ be an adapted c\`adl\`ag process of class
(D), and let $\xi$ be an integrable random variable such that
$L_T\le\xi$. Set
\[
Y_t=\esssup_{\tau\ge t}E(L_\tau\mathbf{1}_{\{\tau<T\}}
+\xi\mathbf{1}_{\{\tau=T\}}|\FF_t).
\]
By Remark \ref{uw3.rev2}, $Y$ is the first component of a solution
$(Y,M,K)$ of RBSDE$^T(\xi,0,L)$. Observe that if $K$ is continuous
then
\begin{equation}
\label{eq3.01} Y_t=E(L_{\tau^*_t}\mathbf{1}_{\{\tau^*_t<T\}}
+\xi\mathbf{1}_{\{\tau^*_t=T\}}|\FF_t),
\end{equation}
where
\[
\tau^*_t=\inf\{s\ge t:Y_s=L_s\}\wedge T.
\]
Indeed, by the definition of $\tau^*_t$ and the definition of a
solution of RBSDE$^T(\xi,0,L)$,
\begin{equation}
\label{eq3.02} Y_t=E\Big((\int_t^{\tau^*_t}\, dK_r
+L_{\tau^*_t}\mathbf{1}_{\{\tau^*_t<T\}}+\xi\mathbf{1}_{\{\tau^*_t=T\}}|\FF_t\Big).
\end{equation}
and, since $K$ is continuous,
\[
\int_0^{T\wedge a}(Y_r-L_r)\,dK_r=0,\quad a\ge 0.
\]
This implies that $\int_t^{\tau^*_t} dK_r=0$, which when combined
with (\ref{eq3.02}) yields (\ref{eq3.01}).
\end{uw}
In the rest of this section we assume that
\begin{equation}
\label{eq3.9}
H^j(t,y)=\max_{i\in A_j}(-c_{j,i}(t)+y^i),
\end{equation}
where $c_{j,i}$ are  continuous adapted process such that for some
constant $c>0$,
\[
c_{j,i}(t)\ge c,\quad i\in A_j,\quad t\in [0,T\wedge a],\quad a\ge
0,\quad j=1,\dots,N.
\]

\begin{uw}
\label{uw3.43} Assume that the underlying filtration $\mathbb{F}$
is quasi-left continuous, $H$ is of the form (\ref{eq3.9}) and $V$
is continuous. Then (A9) is satisfied. Indeed, in this case
$-\Delta K^j_\tau=\Delta Y^j_\tau$ for every predictable stopping
time $\tau$. Therefore repeating step by step the proof of
\cite[Proposition 2]{DHP} we get the desired result.
\end{uw}

\begin{tw}
\label{tw4.4} Assume that $f$ does not depend on $y$ and $H^j$ are
of the form {{\rm(\ref{eq3.9})}}. If
$E(\int_0^T\,d|V|_r+\int_0^T|f(r)|\,dr)<\infty$ then there exists
a solution $(Y,M,K)$ of {{\rm(\ref{eq1.3})}}. Moreover, if $K$ is
continuous, then
\begin{equation}
\label{eq3.10} Y^j_t=\esssup_{\mathcal{S}\in\mathbf{A}_t }
E\Big(\int_t^Tf^{w^j_r}(r)\,dr +\int_t^T\,dV^{w^j_r}_r-\sum_{n\ge
1}c_{w^j_{\tau_{n-1}}, w^j_{\tau_{n}}}
(\tau_n)\mathbf{1}_{\{\tau_n<T\}}+\xi^{w_T}|\FF_t\Big),
\end{equation}
and the optimal strategy $\mathcal{S}^*$ for $Y^j$ is given by
\[
\tau^{j,*}_{0}=0,\quad \xi^{j,*}_{0}=j,\quad \tau_k^{j,*}
=\inf\{t\ge \tau^{j,*}_{k-1}:
Y^{\xi^{j,*}_{k-1}}_t=H^{\xi^{j,*}_{k-1}}_t\}\wedge T,\quad k\ge1,
\]
\[
\xi^{j,*}_k=\sum_{i=1}^N i\mathbf{1}_{\{H^{\xi^{j,*}_{k-1}}_{\tau_k}
=-c_{\xi^{j,*}_{k-1},i}(\tau_k)+Y^i_{\tau_k}\}},\quad k\ge 1.
\]
\end{tw}
\begin{dow}
The existence part follows from Theorem \ref{tw3.2}, because
assumptions (A1)--(A5), (A8) are clearly satisfied, and (A7) is
satisfied with $(\overline{Y},\overline{M})$ and
$(\underline{Y}\,, \underline{M})$ defined as follows:
\[
\overline{Y}^1=\ldots=\overline{Y}^N,\quad
\overline{M}^1=\ldots=\overline{M}^N,\qquad
\underline{Y}^1=\ldots=\underline{Y}^N,\quad
\underline{M}^1=\ldots=\underline{M}^N,
\]
and $(\overline{Y}^1,\overline{M}^1)$  (resp.
$(\underline{Y}^1,\underline{M}^1)$) is a solution of
BSDE$^T(\xi^1\vee\dots\vee\xi^N, f^1\vee\dots\vee
f^N+dV^1\vee\dots\vee dV^N)$ (resp.
BSDE$^T(\xi^1\wedge\dots\wedge\xi^N, f^1\wedge\dots\wedge
f^N+dV^1\wedge\dots\wedge dV^N)$). Thanks to Remark \ref{uw3.4},
to get (\ref{eq3.10}) it suffices now  to repeat step by step the
proof of \cite[Theorem 1]{DHP}.
\end{dow}

\begin{uw}
\label{uw4.444} As a by-product,  from the above theorem we obtain
the following result:  under the assumptions of Theorem
\ref{tw4.4} there exists at most one solution $(Y,M,K)$ of
(\ref{eq1.3}) such that $K$ is continuous. In particular, by
Remark \ref{uw3.43}, if $\mathbb{F}$ quasi-left continuous and $V$
continuous then there is at most one solution of problem
(\ref{eq1.1}) with data satisfying the assumptions of Theorem
\ref{tw4.4}.
\end{uw}

\section{Systems of elliptic quasi-variational inequalities}
\label{sec5}

In this section $E$ is a locally compact separable metric space,
$m$ is a Radon measure on $E$ such that supp$[m]=E$, and
$(\EE,D[\EE])$ is a regular transient semi-Dirichlet form on
$L^2(E;m)$. By $(L,D(L))$ we denote the generator associated with
($\EE,D[\EE])$ (see \cite[Chapter 1]{Oshima}).

Let us recall that $(\EE,D[\EE])$ is called semi-Dirichlet if
$D[\EE]$ is dense in $L^2(E;m)$ and $\EE$ is a bilinear form on
$D[\EE]\times D[\EE]$ satisfying the conditions $(\EE 1)$--$(\EE
4)$ below:
\begin{enumerate}
\item[$(\EE 1)$]  $\EE$ is lower bounded, i.e.
there exists $\alpha_0\ge 0$ such that
\[
\EE_{\alpha_0}(u,u)\ge 0,\quad u\in D[\EE],
\]
where $\EE_{\alpha_0}(u,v)=\EE(u,v)+\alpha_0 (u,v)$,
\item[$(\EE 2)$] $\EE$ satisfies the sector condition, i.e.
there exists $K>0$ such that
\[
|\EE(u,v)|\le K\EE_{\alpha_0}(u,u)^{1/2}\EE_{\alpha_0}(v,v)^{1/2},
\quad u,v\in D[\EE],
\]
\item[$(\EE 3)$] $\EE$ is closed, i.e. for every $\alpha>\alpha_0$
the space $D[\EE]$ equipped with the inner product
$\EE^{(s)}_{\alpha}(u,v):=\frac12(\EE_\alpha(u,v)+\EE_\alpha(v,u))$
is a Hilbert space,
\item[$(\EE 4)$] $\EE$ has the Markov property, i.e. for every $a\ge
0$,
\[
\EE(u\wedge a,u\wedge a)\le \EE(u\wedge a,u),\quad u\in D[\EE].
\]
\end{enumerate}
Note that $(\EE 4)$  is equivalent to the fact that the semigroup
$\{T_t,\,t\ge 0\}$ associated with $(\EE,D[\EE])$ is sub-Markov
(see \cite[Theorem 1.1.5]{Oshima}). Also recall that $\EE$ is said
to have the dual Markov property if
\begin{enumerate}
\item[$(\EE 5)$] for every $a\ge 0$,
\[
\EE(u\wedge a,u\wedge a)\le \EE(u,u\wedge a),\quad u\in D[\EE].
\]
\end{enumerate}
Condition $(\EE 5)$  is equivalent to the fact that associated
dual semigroup $\{\hat{T}_t,t\ge 0\}$ associated with
$(\EE,D[\EE])$ is  sub-Markov (see \cite[Theorem 1.1.5]{Oshima}).
For the notions of transiency and regularity see \cite[Section
1.2, Section 1.3]{Oshima}.

Let Cap denote the capacity associated with $(\EE,D[\EE])$ (see
\cite[Chapter 2]{Oshima}),  and let $\BX=(\{X_t,t\ge 0\}, \{P_x,
x\in E\})$ be a Hunt process with life time $\zeta$ associated
with $(\EE, D[\EE])$ (see \cite[Chapter 3]{Oshima}). We say that
some property holds quasi-everywhere (q.e. for short) if there is
a set $B\subset E$ such that Cap$(B)=0$ and it holds on the set
$E\setminus B$.  A set $B\subset E$ such that Cap$(B)=0$ is called
exceptional.

Let $\mu$ be a signed measure $E$. By $\mu^+$ (resp. $\mu^{-}$) we
denote its positive (resp. negative) part, and we set
$|\mu|=\mu^{+}+\mu^{-}$. A Borel signed measure $\mu$ on $E$ is
called smooth if $\mu$ charges no exceptional sets and there
exists an increasing sequence $\{F_n\}$ of closed subsets of $E$
such that $|\mu|(F_n)<\infty$ for $n\ge 1$, and for every compact
$K\subset E$.
\[
{\rm{Cap}}(K\setminus F_n )\rightarrow 0.
\]
In the sequel the set of all signed smooth measures on $E$ such
that $\|\mu\|_{TV}:=|\mu|(E)<\infty$ will be denoted by
$\MM_{0,b}$.

It is known (see \cite[Section 4.1]{Oshima}) that there is
one-to-one correspondence (the Revuz duality) between positive
continuous additive functionals (PCAFs for short) of $\BX$ and
positive smooth measures. By $A^\mu$ we denote the unique PCAF of
$\BX$ associated with positive smooth measure $\mu$. For a signed
smooth measure $\mu$ we set $A^\mu=A^{\mu^+}-A^{\mu^-}$. By
$\mathbb{M}$ we denote the set of all smooth measures $\mu$ on $E$
such that
\[
E_x\int_0^\zeta\,dA^{|\mu|}_r<\infty
\]
for q.e. $x\in E$, where $E_x$ denotes the expectation with
respect to $P_x$. For a fixed positive measurable function $f$ and
a positive Borel measure $\mu$ we denote by $f\cdot \mu$ the
measure defined as
\[
(f\cdot\mu)(\eta)=\int_E \eta f\,d\mu,\quad \eta\in\mathcal{B}^+(E).
\]
We write  $f\in \mathbb{M}$ if $f\cdot m\in\mathbb{M}$.  By
\cite[Corollary 1.3.6]{Oshima}, if $(\EE,D[\EE])$ has the dual
Markov property then
\begin{equation}
\label{eq3.12} \MM_{0,b}\subset \mathbb{M}.
\end{equation}
By $qL^1(E;m)$ we denote the set of all measurable real functions
$f$ on $E$ such that  $A^{f\cdot m}_t<\infty$ for every $t\ge 0$.
By (\ref{eq3.12}),
\[
L^1(E;m)\subset qL^1(E;m).
\]
Note that in general the form associated with the operator defined
by (\ref{eq1.6}) does not have the dual Markov property.
Nevertheless, for this form (\ref{eq3.12}) holds true.

Recall that a set $U\subset E$ is called quasi-open if for every
$\varepsilon>0$ there exists an open set $U\subset U_\varepsilon
\subset E$ such that Cap$(U_\varepsilon\setminus U)<\varepsilon$.
The family of quasi-open sets induces the quasi-topology on $E$.
We say that a function $u$ on $E$ is quasi-continuous if it is
continuous with respect to the quasi-topology.

\subsection{Existence and approximation of solutions}

For $i,j=1,\dots,N$ let $h_{j,i},\,f^j:E\times\BR^N\rightarrow
\BR$ be measurable functions, $\mu^j$ be  smooth measures  on $E$,
and let $A_j\subset\{1,\dots,j-1,j+1,\dots,N\}$.  We maintain the
notation $f^j(x,y;a)$ introduced at the beginning of Section
\ref{sec2}, and we set
\[
H^j(x,y)=\max_{i\in A_j}h_{j,i}(x,y^i),\qquad H=(H^1,\dots,H^N),
\]
\[
f=(f^1,\dots,f^N),\qquad \mu=(\mu^1,\dots,\mu^N).
\]

We will make the following hypotheses:
\begin{enumerate}[(H1)]
\item $\mu^j\in\mathbb{M}$, $j=1,\dots,N$,
\item for $j=1,\dots,N$ the function $a\mapsto f^j(x,y;a)$
is nonincreasing for  all $x\in E$, $y\in\BR^N$,
\item $f$ is off-diagonal nondecreasing, i.e. for $j=1,\dots,N$ we
have $f^j(x,y)\le f^j(x,\bar{y})$ for all $y,\bar y\in\BR^N$ such
that $y\le\bar{y}$ and $y^j=\bar{y}^j$,
\item  $y\mapsto f(x,y)$ is continuous for every $x\in E$,
\item  $f^j(\cdot,y)\in qL^1(E;m)$ for all $y\in\BR^N$,
$j=1,\dots,N$.
\end{enumerate}

Consider the following system of equations
\begin{equation}
\label{eq4.1} -Lu=f(x,u)+\mu.
\end{equation}

Following \cite{KR:JFA,KR:CM} we adopt the following definition of
a solution of (\ref{eq4.1}).

\begin{df}\label{def3.1}
We say that a measurable function
$u=(u^1,\dots,u^N):E\rightarrow\BR^N$ is a solution of
(\ref{eq4.1}) (PDE$(f+d\mu)$ for short) if
$f^j(\cdot,u)\in\mathbb{M}$, $j=1,\dots,N$, and for q.e. $x\in E$,
\begin{equation}
\label{eq4.12} u^j(x)=E_x\Big(\int_0^\zeta
f^j(X_r,u(X_r))\,dr+\int_0^\zeta\,d A^{\mu^j}_r\Big),\quad
j=1,\dots,N.
\end{equation}
\end{df}

\begin{uw}
A measurable function $u:E\rightarrow\BR^N$ satisfying
(\ref{eq4.12}) may be called a probabilistic solution of
(\ref{eq4.1}). Note that  if $f^j(\cdot,u)\in L^1(E;m)$ and
$\mu^j\in\MM_b$ then $u^j$ is a renormalized solution of
(\ref{eq4.1}) (see \cite{KR:NDEA}).
\end{uw}

\begin{uw}
\label{uw4.1} (i) If $u$ is a solution of (\ref{eq4.1}) in the
sense of Definition \ref{def3.1} then by \cite[Theorem
4.7]{KR:JFA} the pair $(u(X),M)$, where
\begin{equation}
\label{eq4.11} M^j_t=E_x\Big(\int_0^\zeta f^j(X_r,u(X_r))\,dr
+\int_0^\zeta dA^{\mu^j}_r|\FF_t\Big), \quad t\ge0,
\end{equation}
is a solution of BSDE$^\zeta(0,f(X,\cdot)+dA^\mu)$ under the
measure $P_x$ for q.e. $x\in E$ (In fact, $M$ in (\ref{eq4.11}) is
an independent of $x$ version of the right-hand side of equation
(\ref{eq4.11}); such a version always exists, see \cite[Section
A.3]{Fukushima}).
\\
(ii) If $(Y,M)$ is a solution of BSDE$^\zeta(0,f(X,\cdot)+dA^\mu)$
under the  measure $P_x$ for q.e. $x\in E$,
$f^j(\cdot,u)\in\mathbb{M}$, $j=1,\dots,N$, and there exists a
function $u$ such that $u(X)=Y$ under the measure $P_x$ for q.e.
$x\in E$, then $u$ is a solution of (\ref{eq4.1}). This follows
directly from Remark \ref{uw3.0}.
\end{uw}

\begin{df}
We say that a measurable function $u:E\rightarrow\BR^N$ is a
subsolution (resp. supersolution) of (\ref{eq4.1}) if there exists
 positive measures $\beta^1,\dots,\beta^N\in\mathbb{M}$ such that $u$ is a
solution of (\ref{eq4.1}) with $\mu$ replaced by $\mu-\beta$
(resp. $\mu+\beta$), where $\beta=(\beta^1,\dots,\beta^N)$.
\end{df}

\begin{uw}
\label{uw4.2} By Remark \ref{uw4.1}, if $u$ is a subsolution
(resp. supersolution) of (\ref{eq4.1}) then $u(X)$ is the first
component of a subsolution (resp. supersolution) of the equation
BSDE$^\zeta(0,f(X,\cdot)+dA^\mu)$ under the measure $P_x$ for q.e.
$x\in E$.
\end{uw}

\begin{df}
We say that a quasi-continuous function $u$ on $E$ is a solution
of  (\ref{eq1.1}) if there exists   positive measures
$\nu^1,\dots,\mu^N\in\mathbb{M}$ such that $u$ is a solution of
PDE$(f+d\mu+d\nu)$ with $\nu=(\nu^1,\dots,\nu^N)$, and the second
and third condition in (\ref{eq1.1}) are satisfied.
\end{df}

We will also need the following hypotheses:
\begin{enumerate}
\item[(H6)] There exists a subsolution $\underline{u}$ and
a supersolution $\overline{u}$ of (\ref{eq4.1}) such that
\[
\underline{u}\le \overline{u},\quad H(\cdot,\overline{u})
\le \overline{u},\quad \sum_{j=1}^{N}|f^j(\cdot,\overline{u};
\underline{u}^j)|+|f^j(\cdot,\underline{u};\overline{u}^j)|\in\mathbb{M},
\]
\item[(H7)] $H^j$ is continuous  on $E\times\BR^N$
equipped with the product topology consisting of quasi-topology on
$E$ and the Euclidean topology on $\BR^N$ and nondecreasing with respect to $y$.
\end{enumerate}

In the proof of the next theorem we will use some result from
\cite{K:arx} on  solutions of the usual obstacle problem for
single equation and  one quasi-continuous barrier
$h:E\rightarrow\BR$. For the convenience of the reader we recall
below the definition of a solution.

\begin{df}
Let $N=1$. We say that a pair $(u,\nu)$ is a solution of the
obstacle problem for $L$ with lower  barrier $h$ and the
right-hand side $f+d\mu$ (OP$(f+d\mu,h)$ for short) if $u$ is
quasi-continuous, $\nu$ is a positive measure such that $\nu
\in\mathbb{M}$, $u\ge h$ q.e., and
\[
-Lu=f(x,u)+\mu+\nu,\qquad \int_E(u-h)\,d\nu=0.
\]
\end{df}

In the sequel, for $\mu=(\mu^1,\dots,\mu^N)$ we write
$|\mu|=\sum^N_{j=1}|\mu^j|$,
$\|\mu\|_{TV}=\sum^N_{j=1}\|\mu^j\|_{TV}$.
\begin{tw}
\label{tw4.1} Assume {{\rm (H1)--(H7)}}.  Then there exists a
minimal solution of {{\rm(\ref{eq1.1})}} such that
$\underline{u}\le u\le \overline{u}$.
\end{tw}
\begin{dow}
We adopt the notation of  Theorem \ref{tw3.2}. First observe that
the data $f(X,\cdot)$, $H^j(X,\cdot)$, $\xi:=0$, $T:=\zeta$,
$\underline{Y}:=\underline{u}(X)$, $\overline{Y}:=\overline{u}(X)$
satisfy (\ref{eq3.7}) and assumptions (A1)--(A5), (A7), (A8) under
the measure $P_x$ for q.e. $x\in E$ (see Remark \ref{uw4.1}). Set
$u_0=\underline{u}$. By \cite[Theorem 3.2]{K:arx}, for every $n\ge
1$,
\[
u^j_n(X_t)=Y^{n,j}_t,\quad A^{\nu_n}_t=K^{n,j}_t,
\]
where $(u^j_n,\nu^j_n)$  is a solution of
OP$(f^j(\cdot,u_{n-1};\cdot)+\mu^j, H^j(\cdot,u_{n-1}))$. Since
$Y^n\le Y^{n+1}$, we have $u_n\le u_{n+1}$ q.e. Therefore putting
$u^j=\sup_{n\ge 0}u^j_n$ and $u=(u^1,\dots,u^N)$ we obtain
\[
u(X_t)=Y_t,\quad t\in [0,\zeta\wedge a],\quad a\ge 0.
\]
Set $v_n^j(x)=E_x\int_0^\zeta\,dA^{\nu^j_n}_r$, $x\in E$. It is
clear that $v^j_n$ is a potential. From the proof of  Theorem
\ref{tw3.2} it follows that $\{v_n^j\}$ is convergent. Therefore
$v^j:=\lim_{n\rightarrow \infty}v^j_n$ is an excessive function.
By \cite{Oshima}, $v^j$ is a quasi-continuous. Observe that $v^j$
is dominated by the function
$x\mapsto|\overline{u}(x)|+E_x\int_0^\zeta|f(X_r,u(X_r))|\,dr
+E_x\int_0^\zeta\,dA^{|\mu|}_r$. Of course, the the second and the
third term of this sum are potentials. Moreover, $|\overline{u}|$
is dominated by a potential as a supersolution. Thus $v^j$ is
dominated by a potential, and consequently, $v^j$ is a potential.
Therefore  by \cite[Theorem IV (3.13)]{BG} there exists a PCAF
$A^j$ (it is continuous since $v^j$ is quasi-continuous) such that
$v^j(x)=E_\cdot \int_0^\zeta\,dA^j_r$, $x\in E$. By \cite{Oshima}
there exists a positive measure $\nu^j\in\mathbb{M}$ such that
$A^j=A^{\nu^j}$. It is clear that $A^{\nu^j}=K^j$, which implies
that the pair $(u,\nu)$, where $\nu=(\nu^1,\dots,\nu^N)$, is a
solution of (\ref{eq1.1}).
\end{dow}

\begin{uw}
\label{uw4.4} Assume that $(\EE,D[\EE])$ has the dual Markov
property and hypotheses  (H1), (H6), (H7) hold true with
$\mathbb{M}$ replaced by $\MM_{0,b}$. Assume also that the
measures $\beta^j$ appearing in the definition of a supersolution
$\overline{u}$ belongs to $\MM_{0,b}$.  Let $(u,\nu)$ be a
solution of (\ref{eq1.1}). Then $\nu^j\in\MM_{0,b}$,
$j=1,\dots,N$. Indeed, since $u\le\overline{u}$, we have
\begin{align*}
E_x\int_0^\zeta\,dA^{|\nu|}_r&\le E_x\int_0^\zeta|f(r,u(X_r))|\,dr
+E_x\int_0^\zeta|f(X_r,\overline{u}(X_r))|\,dr\\
&\quad+2E_x\int_0^\zeta dA^{|\mu|}+E_x\int_0^\zeta dA^{|\beta|}_r.
\end{align*}
By \cite[Lemma 2.9]{KR:CM}, the above inequality implies that
\[
\|\nu\|_{TV}\le \|f(\cdot,u)\|_{L^1}+ \|f(\cdot,\overline{u})\|_{L^1}
+2\|\mu\|_{TV}+\|\beta\|_{TV}.
\]
By our assumptions, $\|\mu\|_{TV}+\|\beta\|_{TV}<\infty$.
Furthermore, $ \|f(\cdot,u)\|_{L^1}<\infty$ by (H6), and $
\|f(\cdot,\overline{u})\|_{L^1}<\infty$ by (H7), Theorem
\ref{stw3.0} and \cite[Lemma 2.9]{KR:CM}.
\end{uw}

\begin{uw}
Under the assumptions of Remark \ref{uw4.4} the functions $u^j$,
$j=1,\dots$, have the property that $T_k(u^j)\in D_e[\EE]$ for
$k\ge 0$, where $T_k(y)=\max(\min(y,k),-k)$. This follows from
Remark \ref{uw4.4} and \cite[Proposition 5.9]{KR:JFA}. Therefore
under the assumptions of Remark \ref{uw4.4} the function $u^j$ is
a solutions of the first equation in (\ref{eq1.1}) in the sense of
Stampacchia, or, in different terminology, are solution in the
sense of duality (see \cite[Proposition 5.3]{KR:JFA}).
Equivalently, it is a  renormalized solution of this equation (see
\cite{KR:NDEA}).
\end{uw}

\begin{stw}
\label{stw4.2} Let $N=1$. Assume {{\rm (H1), (H3), (H4), (H5)}}.
Moreover, assume that there exists a real valued measurable
function $v$ on $E$ such that $Lv\in\mathbb{M}$ and
$f(\cdot,v)\in\mathbb{M}$. Then there exists a solution $u$ of
PDE$(f+d\mu)$.
\end{stw}
\begin{dow}
Set $\beta=-Au$. Observe that the data $f(X,\cdot),\, V:=A^\mu,\,
\xi:=0,\, S:=v(X),\, T:=\zeta$ satisfy the assumptions of Theorem
\ref{stw3.0} under the measure $P_x$ for q.e. $x\in E$. From the
proof of Theorem \ref{stw3.0} it follows that there exists a
solution $(Y,M)$ of BSDE$^\zeta(0,f(X,\cdot)+dA^\mu)$, and that
$Y=\tilde{Y}+S$, where $(\tilde{Y},\tilde{M})$ is a solution of
BSDE$^T(0,f_S+dA^\mu-dA^\beta)$. By \cite{KR:JFA} there exists a
solution $\tilde{u}$ of PDE$(f_v+d\mu-d\beta)$ with
$f_v(x,y)=f(x,v(x)+y)$,  and $\tilde{u}(X)=\tilde{Y}$. Hence
$Y=\tilde{u}(X)+v(X)$. It is clear (see Remark \ref{uw3.0}) that
$u:=\tilde{u}+v$ is a solution of PDE$(f+d\mu)$.
\end{dow}

In the next proposition we will need the following hypothesis.
\begin{enumerate}
\item[(H9)] There exists a subsolution $\underline{u}$ and a
supersolution $\overline{u}$ of (\ref{eq4.1}), and a measurable
function $v:E\rightarrow\BR^N$such that $Lv^j\in\mathbb{M}$,
$j=1,\dots,N$, and
\[
\underline{u}\le \overline{u},\quad
\sum_{j=1}^{N}|f^j(\cdot,\overline{u};v^j)|
+|f^j(\cdot,\underline{u};v^j)|\in\mathbb{M}.
\]
\end{enumerate}
\begin{stw}
\label{stw4.3} Assume {{\rm (H1)--(H5), (H9)}}. Then there exists
a minimal  solution $u$ of {{\rm(\ref{eq4.1})}} such that
$\underline{u}\le u\le \overline{u}$.
\end{stw}
\begin{dow}
Observe that the data $f(X,\cdot),\, \underline{Y}:=
\underline{u}(X),\, \overline{Y}:=\overline{u}(X),\, S:=v(X),\,
V:=A^\mu,\, \xi:= 0,\, T:=\zeta$ satisfy the assumptions of
Theorem \ref{tw3.2} under the measure $P_x$ for q.e. $x\in E$. Set
$u_0=\underline{u}$. By Proposition \ref{stw4.2} (see also Remark
\ref{uw4.1}), $Y^{j,n}=u^j_n(X)$, where $u^j_n$ is the solution of
PDE$(f^j(\cdot,u_{n-1};\cdot)+d\mu^j)$. From the proof of Theorem
\ref{tw3.2} it follows that $Y^n\le Y^{n+1}$. Hence $u_n\le
u_{n+1}$ q.e. Set $u=\sup_{n\ge 1} u_n$. It is clear that
$Y=u(X)$. Hence,  by Remark \ref{uw4.1}, $u$ is a minimal solution
of PDE$(f+d\mu)$ such that $\underline{u}\le u\le \overline{u}$.
\end{dow}

\begin{tw}
Assume {{\rm (H1)--(H7)}}. Then there exists a minimal  solution
$u_n$ of the system
\begin{equation}
\label{eq4.5} -Lu_n^j=f^j(\cdot,u_n)+n(u_n^j-H^j(\cdot,u_n))^{-}+\mu^j
\end{equation}
such that $\underline{u}\le u_n\le\overline{u}$. Moreover,
$u_n\nearrow u$ q.e.,  where $u$ is the minimal solution of {{\rm
(\ref{eq1.1})}} such that $\underline{u}\le u\le \overline{u}$.
\end{tw}
\begin{dow}
Observe that $\underline{u}$ is a supersolution of (\ref{eq4.5}),
whereas $\underline{u}$ is a subsolution of (\ref{eq4.5}).
Moreover,
$f^j(x,\overline{u})+n(\overline{u}^j-H^j(x,\overline{u}))^-=f^j(x,\overline{u})\in
\mathbb{M}$ by the definition of a supersolution, and
$f^j(x,\underline{u};\overline{u}^j)
+n(\overline{u}^j-H^j(x,\underline{u}))^-=f^j(x,\underline{u};\overline{u}^j)\in
\mathbb{M}$ by (H6). Therefore  (H9) is satisfied for
(\ref{eq4.5}) with $v:=\overline{u}$. By Remark \ref{uw4.1},
$u_n(X)$ is the first component of the solution of
BSDE$^\zeta(0,f_n(X,\cdot)+dA^\mu)$ with
$f^j_n(t,y)=f^j(x,y)+n(y^j-H^j(x,y))^-$. By the construction (see
Proposition \ref{stw4.3}), it is the minimal solution of
BSDE$^\zeta(0,f_n(X,\cdot)+dA^\mu)$. By Theorem \ref{tw3.3} the
sequence $\{u_n(X)\}$ is nondecreasing and $u_n(X)\nearrow Y$,
where $Y$ is the first component of the solution of (\ref{eq1.3})
under the measure $P_x$ for q.e. $x\in E$. From the proof of
Theorem \ref{tw4.1} it follows that $Y=u(X)$, where $(u,\nu)$ is
the minimal solution of (\ref{eq1.1}) such that $\underline{u}\le
u\le \overline{u}$. Of course, this implies that $u_n\nearrow u$
q.e.
\end{dow}

\subsection{Application to the switching problem}

In the theorem below we keep the notation introduced in Section
\ref{sec4}, and we that
\begin{equation}
\label{eq4.6}
H^j(x,y)=\max_{i\in A_j}(-c_{j,i}(x)+y^i),
\end{equation}
where $c_{j,i}$ are  quasi-continuous functions on $E$ such that
for some constant $c>0$,
\[
c_{j,i}(x)\ge c,\quad x\in E,\quad i\in A_j, \quad j=1,\dots,N.
\]

\begin{tw}
Assume that $f$ does not depend on $y$, $H^j$ are of the form
\mbox{\rm (\ref{eq4.6})}, and $f^j,\mu^j\in\mathbb{M}$,
$j=1,\dots,N$. Then there exists a unique solution $u$ of
\mbox{\rm (\ref{eq4.1})}. Moreover,
\[
u^j(x)=\sup_{\mathcal{S}\in\mathbf{A} }E_x\Big(\int_t^\zeta
f^{w^j_r}(X_r)\,dr+\int_t^\zeta\,dA^{\mu^{w^j_r}}_r-\sum_{n\ge1}
c_{w^j_{\tau_{n-1}},
w^j_{\tau_{n}}}(X_{\tau_n})\mathbf{1}_{\{\tau_n<\zeta\}}\Big)
\]
and
\[
u^j(x)=E_x\Big(\int_t^\zeta
f^{w^{j,*}_r}(X_r)\,dr+\int_t^\zeta\,dA^{\mu^{w^{j,*}_r}}_r-\sum_{n\ge
1} c_{w^{j,*}_{\tau_{n-1}}, w^{j,*}_{\tau_{n}}}
(X_{\tau_n})\mathbf{1}_{\{\tau_n<\zeta\}}\Big),
\]
where
\[
w^{j,*}_t=j\mathbf{1}_{[0,\tau^*_1)}(t) +\sum_{n\ge 1}\xi^{j,*}_n
\mathbf{1}_{[\tau^*_{n},\tau^*_{n+1})}(t)
\]
and
\[
\tau^{j,*}_{0}=0,\quad\tau_k^{j,*} =\inf\{t\ge
\tau^{j,*}_{k-1}:u^{\xi^{j,*}_{k-1}}(X_t)
=H^{\xi^{j,*}_{k-1}}(X_t,u(X_t))\}\wedge \zeta,\quad k\ge1,
\]
\[
\xi^{j,*}_{0}=j,\quad \xi^{j,*}_k=\sum_{i=1}^N i
\mathbf{1}_{\{H^{\xi^{j,*}_{k-1}}(X_{\tau_k})
=-c_{\xi^{j,*}_{k-1},i}(X_{\tau_k})+u^i({X_{\tau_k})}\}},\quad
k\ge1.
\]
\end{tw}
\begin{dow}
We know that $\mathbb{F}$ is quasi-left continuous and $A^\mu$ is
continuous. Therefore the theorem follows from Theorem
\ref{tw4.4}, Remark \ref{uw4.444} and Proposition \ref{stw4.3}.
\end{dow}

\bigskip
\noindent{\bf\large Acknowledgements}
\medskip\\
Research supported by the Polish National Science Centre (grant
no. 2012-07-D-ST1-02107).

\end{document}